\documentclass{amsart}

\usepackage{amsmath} 
\usepackage{amssymb}
\usepackage{mathrsfs}

\newtheorem{theorem}{Theorem}[section] 
\newtheorem{claim}[theorem]{Claim}

\newtheorem{conclusion}[theorem]{Conclusion}

\theoremstyle{definition}
\newtheorem{definition}[theorem]{Definition}
\newtheorem{dotes}[theorem]{Definition/Observation}

\newtheorem{choice}[theorem]{Choice}
\newtheorem{discussion}[theorem]{Discussion}
\newtheorem{fact}[theorem]{Fact}
\newtheorem{observation}[theorem]{Observation} 

\newtheorem{hypothesis}[theorem]{Hypothesis}

\theoremstyle{remark}

\newtheorem{remark}[theorem]{Remark}

\newcommand{\rest}{{\restriction}}
\newcommand{\dom}{{\rm dom}} 
 
\newcommand{\Ord}{{\rm Ord}}

\newcommand{\ZFC}{{\rm ZFC}}

\newcommand{\wk}{{\rm wk}}
\newcommand{\fil}{{\rm fil}}
\newcommand{\fin}{{\rm fin}}
\newcommand{\Levy}{{\rm Levy}}
\newcommand{\MA}{{\rm MA}}

\newcommand{\id}{{\rm id}}

\newcommand{\tr}{{\rm tr}}
\newcommand{\st}{{\rm st}}
\newcommand{\wfst}{{\rm wfst}}
\newcommand{\ic}{{\rm ic}}
\newcommand{\cf}{{\rm cf}}
\newcommand{\deep}{{\rm dp}}
\newcommand{\Sp}{{\rm Sp}}
\newcommand{\Dom}{{\rm Dom}}
\newcommand{\Rang}{{\rm Rang}}

\newcommand{\wilog}{{\rm without loss of generality}}

\newcommand{\then}{{\underline{then}}}
\newcommand{\when}{{\underline{when}}}
\newcommand{\Then}{{\underline{Then}}}

\newcommand{\Iff}{{\underline{iff}}}
\newcommand{\mn}{{\medskip\noindent}}
\newcommand{\sn}{{\smallskip\noindent}}

\newcommand{\cA}{{\mathscr A}}

\newcommand{\bbR}{{\mathbb R}}

\newcommand{\gc}{{\mathfrak c}}
\newcommand{\gB}{{\mathfrak B}}

\newcommand{\bbD}{{\mathbb D}}
\newcommand{\bbN}{{\mathbb N}}

\newcommand{\cH}{{\mathscr H}}

\newcommand{\bbL}{{\mathbb L}}

\newcommand{\bbP}{{\mathbb P}}
\newcommand{\cP}{{\mathscr P}}

\newcommand{\bbQ}{{\mathbb Q}}

\newcommand{\cT}{{\mathscr T}}
 
\newcommand{\cU}{{\mathscr U}}

\newcommand{\cY}{{\mathscr Y}}

\newcount\skewfactor
\def\mathunderaccent#1#2 {\let\theaccent#1\skewfactor#2
\mathpalette\putaccentunder}
\def\putaccentunder#1#2{\oalign{$#1#2$\crcr\hidewidth
\vbox to.2ex{\hbox{$#1\skew\skewfactor\theaccent{}$}\vss}\hidewidth}}
\def\name{\mathunderaccent\tilde-3 }

\newenvironment{PROOF}[2][\proofname.]
   {\begin{proof}[#1]\renewcommand{\qedsymbol}{\textsquare$_{\rm #2}$}}
   {\end{proof}}

\begin{document}

\title {The character spectrum of $\beta(\bbN)$}
\author {Saharon Shelah}

\dedicatory  {Dedicated to Kenneth Kunen}

\address{Einstein Institute of Mathematics\\
Edmond J. Safra Campus, Givat Ram\\
The Hebrew University of Jerusalem\\
Jerusalem, 91904, Israel\\
 and \\
 Department of Mathematics\\
 Hill Center - Busch Campus \\ 
 Rutgers, The State University of New Jersey \\
 110 Frelinghuysen Road \\
 Piscataway, NJ 08854-8019 USA}
\email{shelah@math.huji.ac.il}
\urladdr{http://shelah.logic.at}
\thanks{The author thanks Alice Leonhardt for the beautiful typing.\\
Partially supported by the Binational Science Foundation. 
Publication 915}



\date{July 15, 2011}

\begin{abstract}
We show the consistency of: the set of regular cardinals
which are the character of some ultrafilter on $\bbN$ can be quite
chaotic, in particular can have many gaps.
\end{abstract}

\maketitle
\numberwithin{equation}{section}
\setcounter{section}{-1}

\section {Introduction}

The set of characters of non-principal ultrafilters on $\bbN$, that we
call the character spectrum and denote by $\Sp_\chi$, is naturlly of
interest to topologists and set theorists alike, see Definition
\ref{z2} below.  A natural question is what can this set of cardinals
be?  The first result on $\Sp_\chi$ is Posp\'isil proof that $\gc \in
\Sp_\chi$. 

It is consistent that $\Sp_\chi = \{2^{\aleph_0}\}$, since Martin's Axiom
implies $\Sp_\chi = \{2^{\aleph_0}\}$.  Nevertheless, $\Sp_\chi = \{2^{\aleph_0}\}$ is
not a theorem of $\ZFC$.  Juhas (see \cite{Ju80}) proved the
consistency of the existence of a non-principal ultrafilter $D$ so
that $\chi(D) < 2^{\aleph_0}$.  Kunen (in \cite{Ku72}) mentions that
$\aleph_1 \in \Sp_\chi$ in the side-by-side Sacks model.

Those initial resulsts show that $\chi(D)$ is not a trivial cardinal
invariant.  But we may wonder whether $\Sp_\chi$ is an interesting
set.  For instance, can $\Sp_\chi$ include more than two members?
Does it have to be a convex set?  It is proved in \cite[\S6]{BnSh:642} that
$|\Sp_\chi|$ large is consistent, e.g. $2^{\aleph_0}$ is large and all
regular uncountable $\kappa \le 2^{\aleph_0}$ (or just of uncountable
cofinality) belong to it.  It was asked there: among
regular cardinals is it convex?  Now (proved in \cite{Sh:846}) $\Sp_\chi$
does not have to be convex.  In the model of \cite{Sh:846}, there is a
triple of cardinals $(\mu,\kappa,\lambda)$ such that $\mu < \kappa <
\lambda,\mu,\lambda \in \Sp_\chi$ but $\kappa \notin \Sp_\chi$.  In
the present paper we show that $\Sp_\chi$ may exhibit much more
chaotic behavior.

To be specific, starting from two disjoint sets $\Theta_1$ and
$\Theta_2$ of regular uncountable cardinals we produce a forcing
notion $\bbP$ which forces the following properties:
\mn
\begin{enumerate}
\item[$(a)$]  no cardinal (of $\bold V$) is collapsed in $\bold V^{\bbP}$
\sn
\item[$(b)$]  $2^{\aleph_0}$ is an upper bound for the union of $\Theta_1$ and
$\Theta_2$ 
\sn
\item[$(c)$]  $\Theta_1 \subseteq \Sp_\chi$ whereas $\Theta_2 \cap
\Sp_\chi = \emptyset$.
\end{enumerate}
\mn
The proof requires that each element of $\Theta_2$ be measurable and
that $\theta \in \Theta_1$ satisfy $\theta^{< \theta} = \theta$.  This
means that in the extension all members of $\Theta_2$ are weakly
inaccessible and hence that we also do not know for certain that there
are successor cardinals outside $\Sp_\chi$.

In the last section we show that we can, e.g. specify $\Sp_\chi \cap
\aleph_\omega$, basically at will: if we have infinitely many
measurable cardinals then we can make the intersection be
$\{\aleph_n:n \in u\}$ for any subset of $[1,\omega)$ that has no
large gaps, i.e. for every $n$ at least one of $n$ and $n+1$ belongs to
$u$.  If we assume infinitely many compact cardinals then we can
realize any ground model subset of $[1,\omega)$, e.g. $\Sp_\chi \cap
\aleph_\omega$ can be even $\{\aleph_p:p$ prime$\}$.

Let us try to explain how do we do this.  A purpose of
\cite{BnSh:642} is to create a large $\Sp_\chi$.  It provides a way to
ensure many cardinals are in $\Sp_\chi$.  On the other hand,
\cite{Sh:846} provides a way for guaranteeing a cardinal is not
in $\Sp_\chi$.
Here we try to combine the methods, hence creating a large set with
many prescribed gaps which establishes $\Sp_\chi$ in $\bold V^{\bbP}$.

For adding cardinals we use systems of filters, so we deal with them
and with the ``one step forcing" in \S1; we use such systems indexed, e.g. by
$\kappa$-trees, and in the end force by a suitable product of those trees,
not adding reals.  In this direction we do not need large cardinal
assumptions.  For eliminating cardinals we need, essentially,
measurables in the ground model.  After the forcing with $\bbP$, our
measurable cardinals become weakly inaccessible, and we show that they
do not belong to $\Sp_\chi$.

We emphasize that for adding a cardinal to $\Theta_1 \subseteq
\Sp_\chi$, we have to assume $\theta = \theta^{< \theta}$.  Moreover,
$\Theta_2$ consists (in the ground model) of measurable cardinals
which remain weakly inaccessible (= regular limit) cardinals in $\bold
V^{\bbP}$.  Consequently, in \S2 we do not know for certain that there
are successor cardinals outside $\Sp_\chi$.  As in many other cases,
to deal with ``small, e.g. successor" cardinals we have also to
collapse.

The last section of the paper is devoted to the set $\Sp_\chi \cap
\aleph_\omega$.  Let $u \subseteq \omega$ be any set (e.g., $u =
\{p:p$ is a prime number$\}$). If we assume that there are infinitely many
compact cardinals in the ground model, then we can force $\Sp_\chi
\cap \aleph_\omega = \{\aleph_n:n \in u\}$.  Assuming just the
existence of infinitely many measurable cardinals, we can prove a
similar result with some restrictions on $u$.  We need that $|u \cap
\{n,n+1\}| \ge 1$ for every $n \in \omega$.

We thank the referee and Shimoni Garti for helpful comments.

\noindent
Recall
\begin{definition}
\label{z2}
1) For an ultrafilter $D$ on $\bbN$ let $\chi(D)$, the character of
$D$ be min$\{|\cA|:\cA \subseteq D$ and every member of $D$ include
some member of $\cA\}$.

\noindent
2) The character spectrum of non-principal ultrafilters on $\bbN$
   is $\Sp_\chi := \{\chi(D):D$ a non-principal ultrafilter on $\bbN\}$. 
\end{definition}
\newpage

\section {Preliminaries} 

This section is devoted to definitions and facts, needed for proving
the main results of the paper.  We present filter systems $\bar D =
\langle D_t:t \in I\rangle$ and we deal with the one step forcing
$\bbQ_{\bar D}$ where $\bar D = \langle D_\eta:\eta \in {}^{\omega
  >}\omega \rangle,D_\eta$ a filter on $\bbN$ containing the co-finite
subsets of $\bbN$; when $P_1 * \bbQ_{\name{\bar D}_1} \lessdot \bbP_2
* \bbQ_{\name{\bar D}_2}$, and with frames $\bold d = (\bar D_{\bold
  d},F_{\bold d})$ for analyzing $\bbQ_{\bold d}$-name of $\name A$ of
subsets of $\bbN$ modulo the filter on $\bbN$ which $F_{\bold d}$
generated, in particular, a derived $\bbQ_{\bar D}$-name of an ideal
$\name\id_{\bold d}$.

\begin{definition}
\label{z4}
For forcing notion $\bbP_1,\bbP_2$ (i.e. quasi orders).

\noindent
1) $\bbP_1 \subseteq \bbP_2$ \Iff \, $p \in \bbP_1 \Rightarrow p \in
   \bbP_2$ and for every $p,q \in \bbP_1$ we have $\bbP_1 \models ``p \le q"$
   iff $\bbP_2 \models ``p \le q"$.

\noindent
2) $\bbP_1 \subseteq_{\ic} \bbP_2$ \Iff \, $\bbP_1 \subseteq \bbP_2$
   and for every $p,q \in \bbP_1$ we have $p,q$ are compatible in $\bbP_1$
   iff $p,q$ are compatible in $\bbP_2$.

\noindent
3) $\bbP_1 \lessdot \bbP_2$ \Iff \,
\mn
\begin{enumerate}
\item[$\boxplus_1$]   $\bbP_1 \subseteq \bbP_2$ and
   every maximal antichain of $\bbP_1$ is a maixmal antichain of
$\bbP_2$,
\end{enumerate}
\mn
equivalently
\mn
\begin{enumerate} 
\item[$\boxplus_2$]   $\bbP_1 \subseteq_{\ic} \bbP_2$ and 
for every $p_2 \in \bbP_2$ for some $p_1 
\in \bbP_1$ we have $p_1 \le_{\bbP_1} p \Rightarrow$ ($p_2,p$ are
   compatible in $\bbP_2$).
\end{enumerate}
\end{definition}

\begin{dotes}
\label{cn.1}  
1) For ${\cA} \subseteq {\cP}(\bbN)$ let fil$({\cA}) = \{B \subseteq
\omega:\bigcap\limits_{\ell < n} A_\ell \subseteq^* B$ for some 
$n < \omega$ and $A_0,\dotsc,A_{n-1} \in {\cA}\}$; so if $\cA$ is
empty then $\fil(\cA)$ is the filter of co-finite sets.  We may forget
to distinguish between $\cA$ and $\fil(\cA)$.

\noindent
2) fil$({\cA})$ is a filter on $\bbN$ extending the filter of
co-bounded subsets of $\bbN$ but possibly fil$({\cA}) = 
{\cP}(\bbN)$, equivalently $\emptyset \in \fil({\cA})$.

\noindent
3) For a filter $D$ on $X$ let $D^+ = \{Y \subseteq X:Y \ne \emptyset$
   mod $D\}$.
\end{dotes}

\begin{definition}
\label{cn.7}   
Let $I$ be a partial order or just a quasi order.

\noindent
1) We say $\bar D$ is an $I$-filter system \when \,:
\mn
\begin{enumerate}
\item[$(a)$]  $\bar D = \langle D_t:t \in I\rangle$
\sn 
\item[$(b)$]    $D_t \subseteq \cP(\bbN)$ but 
$\emptyset \notin \text{ fil}(D_t)$
\sn
\item[$(c)$]   if $s \le_I t$ then fil$(D_s) \subseteq 
\text{ fil}(D_t)$.
\end{enumerate}
\mn
2) We say $\bar D$ is an ultra $I$-filter system when in addition:
\mn
\begin{enumerate} 
\item[$(d)$]  if $s \in I,A \subseteq \bbN$ and $A \ne \emptyset$
mod $D_s$ \then \, for some $t$ we have $s \le_I t$ and 
$A \in \text{ fil}(D_t)$.
\end{enumerate}
\mn
3)  If $\bar D_\ell$ is an $I_\ell$-filter system for $\ell=1,2$ \then
\, we let ($\bar D_\ell = \langle D_{\ell,t}:t \in I_\ell\rangle$ and):
\mn
\begin{enumerate}
\item[$(a)$]    $\bar D_1 \le \bar D_2$ \underline{means} $I_1 \subseteq
I_2$ (as quasi orders, so possibly $I_1=I_2$) and
$s \in I_1 \Rightarrow D_{1,s} \subseteq D_{2,s}$
\sn
\item[$(b)$]   $\bar D_1 \le^* \bar D_2$ means $I_1 \subseteq I_2$ and
$s \in I_1 \Rightarrow \fil(D_{1,s}) \subseteq \fil(D_{2,s})$
\sn
\item[$(c)$]   $\bar D^1 \le^\odot \bar D^2$ means $I_1 \subseteq I_2$ and
$s \in I_1 \Rightarrow \fil(D_{1,s}) = \fil(D_{2,s})$
\sn
\item[$(d)$]   $\bar D^1 =^* \bar D^2$ means $I_1 = I_2$ and
$s \in I_1 \Rightarrow \fil(D_{1,s}) = \fil(D_{2,s})$.
\end{enumerate}
\end{definition}

\begin{observation}
\label{cn.14}  
Let $I$ be a partial order.

\noindent
0) $\le,\le^\odot$ and $\le^*$ quasi order the set of $I$-filter systems
and $\langle \fil(D_t):t \in I\rangle$ is an $I$-filter system for any
$I$-filter system $\bar D$ and $\bar D_1 \le \bar D_2 \Rightarrow \bar
D_1 \le^* \bar D_2$ and $\bar D^1 =^* \bar D^2 \Rightarrow \bar D^1
\le^\odot \bar D^2 \Rightarrow \bar D^1 \le^* \bar D^2$ and $\bar D_1
\le \bar D_2 \le D_1 \Rightarrow \bar D_1 = \bar D_2$ and $\bar D_1
\le^* \bar D_2 \le^* \bar D_1 \Rightarrow \bar D_1 =^* \bar D_2$.

\noindent
1) If $A_s \in [\bbN]^{\aleph_0}$ for each $s \in I$ 
and $A_t \subseteq^* A_s$ for $s
\le_I t$ \then \, there is an $I$-filter system $\bar D$ such that $s \in I
\Rightarrow D_s = \{A_s\}$.

\noindent
2) If $\bar D$ is an $I$-filter system \then \, for some ultra
$I$-filter system $\bar D'$ we have $\bar D \le \bar D'$.

\noindent
3) If $\bar D$ is an $I$-filter system, $s \in I$ and $A \subseteq
\omega$ and $(\forall t)[s \le_I t \Rightarrow A \ne \emptyset$ {\rm mod fil}
$(D_t)]$, \then \, for some $I$-filter system $\bar D'$ we have
$\bar D \le \bar D'$ and $A \in D'_s$.

\noindent
4) If $\langle \bar D_\alpha:\alpha < \delta\rangle$ is an $\le$-increasing
sequence of $I$-filter systems \then \, some $I$-filter system $\bar
D_\delta$ is an upper bound of the sequence; in fact, one can use the limit,
i.e. $D_{\delta,s} = \cup\{D_{\alpha,s}:\alpha < \delta\}$; 
similarly for $\le^*$-increasing.

\noindent
5) If $\bar D$ is an $I$-filter system and $\bar D' = \langle\text{\rm
fil}(D_t):t \in I\rangle$ \then \, $\bar D \le \bar D'$.

\noindent
6) If $\bar D$ is an $I$-filter system and each $D_t$ is an
ultrafilter on $\omega$ \then \, $\bar D$ is an ultra $I$-filter system and
necessarily $s \le_I t \Rightarrow D_s = D_t$.

\noindent
7) If $\bar D_1$ is an ultra $I$-filter system and $\bar D_2$ is an
$I$-filter system such that $\bar D_1 \le^* \bar D_2$ \then \, $\bar D_1
\le^\odot \bar D_2$.

\noindent
8) Assume $\bbP_1 \lessdot \bbP_2$ and $\Vdash_{\bbP_1} ``\name{\bar
   D}_\ell$ is an $I$-filter system" for $\ell=1,2$.  If
   $\Vdash_{\bbP_1} ``\name{\bar D}_1 \le \name{\bar D}_2"$ \then \,
   $\Vdash_{\bbP_2} ``\name{\bar D}_1 \le \name{\bar D}_2"$; also if
   $\Vdash_{\bbP_1} ``\name{\bar D}_1 \le^* \name{\bar D}_2$ \then \,
   $\Vdash_{\bbP_2} ``\name{\bar D}_1 \le^* \name{\bar D} _2"$. 

\noindent
9) If $\bbP_1 \lessdot \bbP_2$ and $\Vdash_{\bbP_1} ``\name{\bar
   D}_\ell$ is an $\name I_\ell$-filter system" for $\ell=1,2$ and
$\Vdash_{\bbP_1} ``\name{\bar D}_1$ is ultra and 
$\name{\bar D}_1 \le^* \name{\bar D}_2"$ \then \,
$\Vdash_{\bbP_2} ``\name D_{1,t} \subseteq \name D_{2,t}$ and (fil$(\name
D_{1,t})^+)^{\bold V[\bbP_1]} \subseteq \text{ fil}(\name D_{2,t})^+"$.
\end{observation}

\begin{PROOF}{\ref{cn.14}}
0) Easy.

\noindent
1) Check.

\noindent
2) Use parts (3),(4), easy, but we elaborate.  We try to 
choose $\bar D_\alpha$ by induction on $\alpha < (2^{\aleph_0} +
|I|)^+$ such that $\bar D_\alpha$ is an $I$-filter system, $\beta <
\alpha \Rightarrow \bar D_\beta \le \bar D_\alpha$ and for each
$\alpha = \beta +1$ for some $t,D_{\alpha,t} \ne D_{\beta,t}$.  For
$\alpha=0$ let $\bar D_\alpha = \bar D$, for $\alpha$ limit use part
(4) and for $\alpha = \beta +1$ if $\bar D_\beta$ is not ultra, use
part (3).  By cardinality consideration for some $\beta,\bar D_\beta$
is defined but we cannot define $\bar D_{\beta +1}$ so necessarily
$\bar D_\beta$ is ultra as required.

\noindent
3)-9) Easy, too. 
\end{PROOF}

\begin{claim}
\label{cn.17}  
1) Assume the quasi-order $I$ as a forcing notion adds no new reals.
An $I$-filter system $\bar D$ is ultra \Iff \, $\Vdash_I
``\cup\{\fil(D_t):t \in \bold G_I\}$ is an ultrafilter on $\omega$".

\noindent
2) Assume the quasi-order $I$ as a forcing notion adds no new
$\omega_1$-sequences of ordinals and $\bbP$ is a {\rm c.c.c.} 
forcing notion, (or just $I$ is $\aleph_1$-complete or just
$\Vdash_{\bbP}$ ``forcing with $I$ add no new
real").  If $\Vdash_{\bbP} ``\langle \name D_t:t \in I\rangle$ is an
$I$-filter system" \then \, $\Vdash_{\bbP} \Vdash_I 
``\cup\{\fil(\name D_t):t \in \name G_I\}$ is an ultrafilter on $\bbN"$ iff
$\Vdash_{\bbP} ``\langle \name D_t:t \in I\rangle$ is an ultra 
$I$-filter system".
\end{claim}

\begin{PROOF}{\ref{cn.17}}
Easy.
\end{PROOF}

\begin{discussion}
\label{cn.19}
An $I$-filter system $\bar D$ may be ``degenerated", i.e. $D_t = D$
   is an ultrafilter, the same for every $t \in I$.  But in this case adding a
generic set to $I$ will not add naturally a new ultrafilter, which is our aim
   here.
\end{discussion}

\begin{definition}
\label{c20}
1) For $\bar D = \langle D_\eta:\eta \in {}^{\omega >}\omega \rangle$,
 each $D_\eta$ a filter on $\bbN$ let $\bbQ_{\bar D}$ be

\begin{equation*}
\begin{array}{clcr}
\bigl\{T:&T \subseteq {}^{\omega >} \omega \text{ is closed under
initial segments, and for some} \\
   &\text{tr}(T) \in {}^{\omega >} \omega,
\text{ the trunk of } T, \text{ we have}: \\
   &(i) \quad \ell \le \ell g(\text{tr}(T)) \Rightarrow T \cap {}^\ell
\omega = \{\text{tr}(T) \restriction \ell\} \\
   &(ii) \quad \text{tr}(T) \trianglelefteq \eta \in {}^{\omega >} \omega
\Rightarrow \{n:\eta \char 94 \langle n \rangle \in T\} \in \name
D_\eta \bigr\}
\end{array}
\end{equation*}

\mn
ordered by inverse inclusion.

\noindent
2) For $p \in \bbQ_{\bar D}$ let wfst$(p,\bar D)$ be the set of pairs
   $(S,\zeta)$ such that:
\mn
\begin{enumerate}
\item[$(a)$]   $(\alpha) \quad S \subseteq \{\eta \in p:\tr(p)
\trianglelefteq \eta \in p\}$
\sn
\item[${{}}$]  $(\beta) \quad \tr(p) \in S$
\sn
\item[${{}}$]  $(\gamma) \quad \tr(p) \trianglelefteq \nu
\triangleleft \eta \in S \Rightarrow \nu \in S$
\sn
\item[$(b)$]  $(\alpha) \quad \zeta$ is a function from $S$ into $\omega_1$
\sn
\item[${{}}$]  $(\beta) \quad$ if $\nu \triangleleft \eta$ are
from $S$ then $\zeta(\nu) > \zeta(\eta)$
\sn
\item[${{}}$]  $(\gamma) \quad$ if $\eta \in S$ and
$\zeta(\eta) > 0$ \then \, $\{k:\eta \char 94 \langle k \rangle 
\in S\} \ne \emptyset \mod D_\eta$.
\end{enumerate}
\mn
3) If $p \in \bbQ_{\bar D}$ and $\nu \in p$ then we let $p^{[\nu]} =
\{\rho \in p:\rho \trianglelefteq \nu$ or $\nu \trianglelefteq
\rho\}$.

\noindent
4) If $\bar D = \langle D_\eta:\eta \in {}^{\omega
>}\omega\rangle,D_\eta = D$ for $\eta \in {}^{\omega >}\omega$ then
let $\bbQ_D = \bbQ_{\bar D}$ and $\wfst(p,D) = \wfst(p,\bar D)$; we may write
$\eta$ instead of $p$ when this holds for some $p \in \bbQ_{\bar D}$ with
$\tr(p) = \eta;\wfst$ stands for well founded sub-tree. 
\end{definition}

\begin{claim}
\label{c30}  
Assume $\eta^* \in {}^{\omega >} \omega,D_\eta$ is a filter on $\bbN$
for $\eta \in {}^{\omega >}\omega$ and 
$\cY$ is a subset of $\Lambda = \Lambda_{\eta^*}
= \{\eta:\eta^* \trianglelefteq \eta \in {}^{\omega >} \omega\}$.  
\Then \, exactly one of the following clauses holds:
\mn
\begin{enumerate}
\item[$(a)$]  there is $q \in \bbQ_{\bar D}$ such that
\begin{enumerate}
\item[$(\alpha)$]   $\eta^* = \tr(q)$
\sn 
\item[$(\beta)$]  ${\cY} \cap q = \emptyset$, equivalently 
$q^+ = q \backslash \{\tr(q(t)) \restriction \ell:\ell < 
\ell g(\tr(q))\})$ is disjoint to ${\cY}$
\end{enumerate}
\item[$(b)$]  there is a function $\zeta$ such that
$(\Dom(\zeta),\zeta) \in \wfst(\eta^*,\bar D)$ and $\max(\Dom(\zeta))
\subseteq \cY$; that is:
\sn
\begin{enumerate}
\item[$(\alpha)$]   $\Dom(\zeta)$ is a set $\Xi$ satisfying
\sn
\item[${{}}$]    $(i) \, \Xi \subseteq \{\eta:\eta^* \trianglelefteq 
\eta \in {}^{\omega >} \omega\}$
\sn
\item[${{}}$]   $(ii) \,\,\eta^* \in \Xi$
\sn
\item[${{}}$]   $(iii) \,\,$ if $\eta \in \Xi$ and
$\eta^* \trianglelefteq \nu \trianglelefteq \eta$ then $\nu \in \Xi$
\sn
\item[$(\beta)$]   $(i) \quad \Rang(\zeta) \subseteq \omega_1$
\sn
\item[${{}}$]  $(ii) \,\, \eta^* \trianglelefteq \nu \triangleleft 
\eta \in \Xi \Rightarrow \zeta(\eta) < \zeta(\nu)$ 
\sn
\item[$(\gamma)$]   for every $\eta \in \Xi$ at 
least one of the following holds:
\sn
\item[${{}}$]   $(i) \quad \eta \in {\cY}$
\sn
\item[${{}}$]  $(ii) \,\,$ the set
$\{n:\eta \char 94 \langle n \rangle \in \Xi\}$ belongs to $D^+_\eta$.
\end{enumerate}
\end{enumerate}
\end{claim}

\begin{PROOF}{\ref{c30}}
Similar to \cite[4.7]{Sh:700} or better \cite[5.4]{Sh:707}.

In full, recall $\Lambda = \{\eta:\eta^* \trianglelefteq \eta \in
{}^{\omega >}\omega\}$.  We define when $\deep(\eta) \ge \zeta$ for
$\eta \in \Lambda$ by induction on the ordinal $\zeta$:
\mn
\begin{enumerate}
\item[$\boxplus$]  $\bullet \quad$ \underline{$\zeta =0$}:  always
\sn
\item[${{}}$]   $\bullet \quad$ \underline{$\zeta$ a limit ordinal}:
$\deep(\eta) \ge \zeta$ iff rk$(\eta) \ge \xi$ for every $\xi < \zeta$
\sn
\item[${{}}$]   $\bullet \quad$ \underline{$\zeta = \xi+1$}:  $\deep(\eta)
\ge \zeta$ \Iff \, both of the following occurs:
\sn
\begin{enumerate}
\item[${{}}$]   $(i) \quad \eta \notin \cY$ 
\sn
\item[${{}}$]   $(ii) \,\,$ the following set belongs 
to $D^+_\eta:\{n:\deep(\eta \char 94 \langle n \rangle) \ge \xi\}$.
\end{enumerate}
\end{enumerate}
\mn
We define $\deep(\eta) \in \Ord \cup \{\infty\}$ such that $\xi = 
\deep(\eta)$
\Iff \, $(\forall \zeta \in \Ord)[\deep(\eta) 
\ge \zeta$ iff $\zeta \le \xi]$.

Easily
\mn
\begin{enumerate}
\item[$\boxplus$]  for every $\eta \in \Lambda,\deep(\eta) 
\in \omega_1 \cup \{\infty\}$.
\end{enumerate}
\bigskip

\noindent
\underline{Case 1}:  $\deep(\eta^*) = \infty$.

For each  $\eta \in \Lambda$ such that $\deep(\eta) = \infty$ clearly
there is $A_\eta \in D_\eta$ such that $n \in A_\eta \Rightarrow
\deep(\eta \char 94 \langle n \rangle) = \infty$.  Let $q$ be

\[
\{\nu \in p:\nu \trianglelefteq \eta^* \text{ or } \eta^*
\triangleleft \nu \text{ and if } \eta^* \trianglelefteq \rho
\triangleleft \nu \text{ then } \nu(\ell g(\rho)) \in A_\rho\}.
\]

\mn
Clearly $q$ is as required in clause (a) of \ref{c30}.
\bigskip

\noindent
\underline{Case 2}:  $\deep(\eta^*) < \infty$.

We define

\begin{equation*}
\begin{array}{clcr}
\Xi = \{\nu:& \eta^* \trianglelefteq \nu \text{ and if } k \in [\ell
g(\eta^*),\ell g(\nu)) \\
  &\text{ then } \nu \rest k \notin \cY \text{ and } \deep(\nu \rest
k) > \deep(\nu \rest (k+1))\}.
\end{array}
\end{equation*}

\mn
We define $\zeta:\Xi \rightarrow \omega_1$ by $\zeta(\eta) =
\deep(\eta)$.

Now check.
\end{PROOF}

\begin{claim}
\label{cn.28}
$\bbP_1 * \bbQ_{\name D_1} \lessdot \bbP_2 * \bbQ_{\name D_2}$ \when \,:
\mn
\begin{enumerate}
\item[$(a)$]   $\bbP_1 \lessdot \bbP_2$ and $\name{\bar D}_\ell = \langle
\name D_{\ell,\eta}:\eta \in {}^{\omega >}\omega\rangle$ for $\ell=1,2$
\sn
\item[$(b)$]  $\name D_{1,\eta}$ is a $\bbP_1$-name of a filter on $\bbN$
\sn
\item[$(c)$]   $\name D_{2,\eta}$ is a $\bbP_2$-name of a filter on $\bbN$
\sn
\item[$(d)$]   $\Vdash_{\bbP_2}$ ``$\name D_{1,\eta} \subseteq 
\name D_{2,\eta}$ and moreover $(\fil(\name D_{1,\eta})^+)^{\bold V[\bbP_1]} 
\subseteq \text{\rm fil}(\name D_{2,\eta})^+$,
i.e. for every $A \in {\cP}(\bbN)^{\bold V[\bbP_1]}$ we have 
$A \in \fil(\name D_{1,\eta}) \Leftrightarrow A \in \fil(\name D_{2,\eta})"$.
\end{enumerate}
\end{claim}

\begin{PROOF}{\ref{cn.28}}
Like \cite[\S4]{Sh:700} more \cite[\S5]{Sh:707} but we elaborate.

Without loss of generality $\emptyset \in \bbP_1$ and $\emptyset
\le_{\bbP_2} p$ for every $p \in \bbP_2$.
Clearly $\bbP_1 * \bbQ_{\name{\bar D}_1} \subseteq \bbP_2 *
\bbQ_{\name{\bar D}_2}$ by clause (d) of the assumption and 
moreover $\bbP_1 \lessdot \bbP_2 \lessdot \bbP_2 * 
\bbQ_{\name{\bar D}_2}$ recalling Definition \ref{z4}(1),(2).  
Now we can force by $\bbP_1$ 
so \wilog \, it is trivial, hence we have to prove that 
$\bbQ_{\bar D_1} \lessdot \bbP_2 * \bbQ_{\name{\bar D}_2}$ identifying
$q \in \bbQ_{\bar D_1}$ with $(\emptyset,q) \in \bbP_2 *
\name{\bbQ}_{\name{\bar D}_2}$.  By clause (d) of the assumption, this
identification is well defined and $\bbQ_{\bar D_1} \subseteq_{\ic}
\bbP_2 * \bbQ_{\name{\bar D}_2}$ because for $p_1,p_2 \in \bbQ_{\bar
  D_1},p_1,p_2$ are compatible iff $(\tr(p_1) \in p_2) \vee (\tr(p_2)
\in p_1)$.  It suffices to verify \ref{z4}(3), requirement $\boxplus_2$.
So let $(p_2,\name q_2) \in \bbP_2 * \bbQ_{\name D_2}$; \wilog
   \, for some $\eta^*$ from $\bold V$ we have $p_2 
\Vdash ``\eta^* = \text{tr}(\name q_2)"$, so 
$\eta^* \in {}^{\omega >}\omega$ and of course:
\mn
\begin{enumerate}
\item[$(*)_1$]  $\Vdash_{\bbP_2} ``\name q_2 \in \bbQ_{\bar D_2}"$.
\end{enumerate}
\mn
By \ref{z4}(3), it suffices to find $q \in \bbQ_{\bar D_1}$ such that
\mn
\begin{enumerate}
\item[$(*)_2$]  $q \le q' \in \bbQ_{\bar D_1} \Rightarrow (p_2,\name q_2),q'$
are compatible; that is, $(p_2,\name q_2),(\emptyset,q')$ are
compatible in $\bbP_2 * \bbQ_{\name{\bar D}_2}$.
\end{enumerate}
\mn
Now we shall apply Claim \ref{c30} in $\bold V$ 
with $\eta^*,\bar D_1$ here standing for $\eta^*,\bar D$ there.
Still $\cY$ is missing, so let

\begin{equation*}
\begin{array}{clcr}
\cY = \{\nu: &\eta^* \trianglelefteq \nu \in {}^{\omega >}\omega \\
  &\text{ and there is } r \in \bbQ_{\bar D_1} \text{ such that } 
\nu = \tr(r) \text{ and} \\
  &(\emptyset,r),(p_2,\name q_2) \text{ are incompatible in } 
\bbP_2 * \bbQ_{\name{\bar D}_2} \\
  &\text{ equivalently } p_2 \Vdash_{\bbP_2} ``\name q_2,r \text{ are
incompatible in } \name{\bbQ}_{\name{\bar D}_2}"\}.
\end{array}
\end{equation*}

\mn
By Claim \ref{c30} below we get clause (a) or clause (b) there.
\bigskip

\noindent
\underline{Case 1}:  Clause (a) holds, say 
as witnessed by $q \in \bbQ_{\bar D_1}$.

We shall prove that in this case $q$ is as required, i.e. $q
\in \bbQ_{\bar D_1}$ and $[q \le_{\bbQ_{\bar D_1}} r \in 
\bbQ_{\bar D_1} \Rightarrow 
(p_2,\name q_2) \in \bbP_2 * \bbQ_{\name{\bar D}_2}$ and $r$ are
compatible (in $\bbP_2 * \bbQ_{\name{\bar D}_2})]$.

Why?  Let $\nu = \tr(r)$.  Clearly
$(\eta^* \trianglelefteq \nu \in q)$ hence
by the choice of $q$, i.e. \ref{c30}$(a)(\beta)$ we have $\nu
\notin \cY$ so $r$ cannot witness $``\nu \in \cY"$ hence
$r,(p_2,\name q_2)$ are compatible in $\bbP_2 * \bbQ_{\name{\bar
D}_2}$ as required.
\bigskip

\noindent
\underline{Case 2}:  Clause (b) holds as witnessed by the function
$\zeta$.

By the definition of $\cY$, in $\bold V$, we can choose $\bar q$ such that:
\mn
\begin{enumerate}
\item[$\boxplus$]  $(a) \quad \bar q = \langle q_\nu:\nu \in
\cY\rangle$
\sn
\item[${{}}$]  $(b) \quad q_\nu \in \bbQ_{\bar D_1}$ and $\tr(q_\nu) =
\nu$
\sn
\item[${{}}$]  $(c) \quad q_\nu$ witness $\nu \in \cY$, i.e. $p_2
\Vdash ``q_\nu,\name q_2$ are incompatible in $\bbQ_{\name{\bar D}_2}"$.
\end{enumerate}
\mn
We define a $\bbP_2$-name $\name q_*$ as follows:

\begin{equation*}
\begin{array}{clcr}
\name q_* = \{\nu:&\nu \trianglelefteq \eta^* \text{ \underline{or} }
\eta^* \triangleleft \nu \in \name q_2 \text{ and if } 
\ell g(\eta^*) \le k < \ell g(\nu) \\
 &\text{ and } \nu \rest k \in \cY \text{ then } \nu \in q_{\nu \rest k}, 
\text{ hence} \\
  &k \le \ell \le \ell g(\nu) \Rightarrow \nu \rest \ell \in
q_{\nu \rest k}\}.
\end{array}
\end{equation*}

\mn
Clearly $\Vdash_{\bbP_2} ``\name q_* \in \bbQ_{\name{\bar D}_2}$ and
$\tr(\name q_*) = \eta^*$ and $\bbQ_{\name{\bar D}_2} \models ``\name
q_2 \le \name q_*""$.
\mn
\begin{enumerate}
\item[$(*)_3$]  if $\nu \in \cY$ then $\eta^* \trianglelefteq \nu$ and
$p_2 \Vdash_{\bbP_2} ``\neg(\nu \in \name q_*)"$.
\end{enumerate}
\mn
[Why?  Otherwise there is $p_3 \in \bbP_2$ such that $p_2 \le p_3$ and
$p_3 \Vdash_{\bbP_2} ``\eta^* \trianglelefteq 
\nu \in \name q_*"$, as $\tr(\name
q_*)$ is forced to be $\eta^*$ and $\tr(q_\nu) = \nu$, necessarily
$p_3 \Vdash_{\bbP_2} ``q_\nu,\name q_*$ are compatible".  But $p_2
\Vdash_{\bbP_2} ``\name q_2 \le \name q_*"$, we get a contradiction to
the choice of $q_\nu$.]

Now we know that $\eta_* \in \Dom(\zeta)$ and 
$\Vdash ``\eta^* \in \name q_*"$ hence $S :=
\{\nu:\nu \in \Dom(\zeta)$ hence $\eta^* \trianglelefteq \nu$ and $p_2
\nVdash ``\nu \notin \name q_*"\}$ is not empty.  So as $S \subseteq
\Dom(\zeta)$ the set $\cU = \{\zeta(\nu):\nu \in S\}$ is not empty, and by
the choice of the function $\zeta$ we have $\cU \subseteq \omega_1$, 
hence there is a minimal $\gamma \in
\cU$ and let $\nu \in \Dom(\zeta)$ be such that $\zeta(\nu) = \gamma$.  By
the definition, if $\gamma=0$ then by clauses $(\gamma)$ and $(\beta)$
of \ref{c30}(b), i.e. the choice of $\zeta(-)$ we have
 $\nu \in \cY$ and, of course, $\nu \in S$.  By $(*)_3$, $p_2 \Vdash_{\bbP_2}
``\neg(\nu \in \name q_*)"$ we get easy contradiction to $\nu \in S$,
hence we can assume
$\gamma > 0$.  By the definition of $S$ there is $p_* \in \bbP_2$ such
that $\bbP_2 \models ``p_2 \le p_*"$ and $p_* \Vdash_{\bbP_2} ``\nu
\in \name q_*$ hence $\in \name q_2"$ and, of course, $\nu \in S$.  
By the choice of the function $\zeta$, in $\bold V$ we have 
$A := \{n:\nu \char 94 \langle n \rangle \in
\Dom(\zeta)\} \ne \emptyset$ mod $D_{1,\nu}$, hence by clause (d) of 
the assumption
of the claim $\Vdash_{\bbP_2} ``A \ne \emptyset$ mod $\name D_{2,\nu}"$
and, of course, $p_* \Vdash_{\bbP_2} ``\{n:\nu \char 94 \langle n
\rangle \in \name q_*\} \in \name D_{2,\nu}"$.  Together $p_*
\Vdash_{\bbP_2}$ ``there is $n$ such that $\nu \char 94 \langle n
\rangle \in \name q_* \cap \Dom(\zeta)"$, so let $n_*$ and $p_{**} \in
\bbP_2$ be such that $\bbP_2 \models ``p_* \le p_{**}"$ and $p_{**}
\Vdash_{\bbP_2} ``\nu \char 94 \langle n_* \rangle \in \name q_* \cap
\Dom(\zeta)"$.

So $\zeta(\nu \char 94 \langle n_* \rangle)$ is well defined,
i.e. $\nu \char 94 \langle n_* \rangle$ belongs
to $\Dom(\zeta)$ hence $\zeta(\nu \char 94 \langle n_* \rangle) < \zeta(\nu) =
\gamma$ and easily $\nu \char 94 \langle n_* \rangle \in S$ and 
$\zeta(\nu \char 94 \langle n_* \rangle) \in \cU$, so
we get a contradiction to the choice of $\gamma$.
\end{PROOF}

\begin{definition}
\label{c31}
1) We say $\bold d = (\bar D,F)$ is a frame \when \,:
\mn
\begin{enumerate}
\item[$(a)$]  $\bar D = \langle D_\eta:\eta \in {}^{\omega >}\omega
\rangle$ and $D_\eta \subseteq [\bbN]^{\aleph_0},\emptyset \notin
\fil(D_\eta)$ for $\eta \in {}^{\omega >}\omega$
\sn
\item[$(b)$]  $F \subseteq [\bbN]^{\aleph_0}$ and $\emptyset \notin
\fil(F)$.
\end{enumerate}
\mn
1A) Above let $\bar D_{\bold d} = \langle D_{\bold d,\eta}:\eta \in
{}^{\omega >}\omega\rangle,D_{\bold d,\eta} =
\fil(D_\eta),F_{\bold d} = \fil(F),\bbQ_{\bold d} = \bbQ_{\bar
D_{\bold d}}$ and if $D_\eta = D$ for $\eta \in {}^{\omega >}\omega$ we may
write $D,D_{\bold d}$ instead of $\bar D,\bar D_{\bold d}$, respectively. 

\noindent
2) We say $\name A$ is a $\bold d$-candidate \when \, ($\bold d$ is a
   frame and):
\mn
\begin{enumerate}
\item[$(c)$]  $\name A$ is a $\bbQ_{\bold d}$-name of a subset of
$\bbN$.
\end{enumerate}
\mn
3) We say $\name A$ is $\bold d$-null when it is a $\bold d$-candidate
and is not $\bold d$-positive, see below.

\noindent
4) We say $\name A$ is $\bold d$-positive \when \, for some $p_* \in
\bbQ_{\bold d}$, for a dense set of $p \ge p_*$ some quadruple 
$(p,A,\bar S,\bar \zeta)$ is a local witness\footnote{An equivalent
  version is when we weaken clause (e) to: if $\eta \in S_n$ and
  $\zeta_n(\eta) = 0$ then there is $q \in \bbQ_{\bold d}$ such that
  $\tr(q) = \eta,p \le q$ and $q \Vdash ``n \in \name A"$, see
  $(*)_{2.2}$ in the proof.  Moreover, we can omit ``$p \le q$"; hence
actually only $\tr(p)$ is important so we may write $\tr(p)$ instead of $p$.}
for $(\name A,\bold d)$ \underline{or} for $(\eta,\name A,\bold d)$ when $\eta
= \tr(p)$ \underline{or} for $(p,\name A,\bold d)$ \underline{or} 
for $\name A$ being $\bold d$-positive, which means:
\mn
\begin{enumerate}
\item[$(a)$]  $p \in \bbQ_{\bold d}$
\sn
\item[$(b)$]  $A \in F^+_{\bold d}$
\sn
\item[$(c)$]  $\bar S = \langle S_n:n \in A\rangle$ and
$\bar\zeta = \langle \zeta_n:n \in A\rangle$
\sn
\item[$(d)$]  $(S_n,\zeta_n) \in \wfst(p,\bar D)$ for $n \in A$ recalling
Definition \ref{c20}(2)
\sn
\item[$(e)$]  if $\eta \in S_n$ and
$\zeta_n(\eta)=0$ then $p^{[\eta]} \Vdash ``n \in \name A"$.
\end{enumerate}
\end{definition}

\begin{definition}
\label{c33}
1) For a frame $\bold d = (\bar D,F)$ let $\name\id_{\bold d} =
\name\id(\bold d) = 
\{\name A \subseteq \bbN:\name A$ is a $\bbQ_{\bold d}$-name 
which is $\bold d$-null$\}$.

\noindent
2) If $\Vdash_{\bbP} ``\name{\bold d}$ is a frame" then 
$\name \id_{\name{\bold d}}[\bbP]$ is the $\bbP * \bbQ_{\name{\bold d}}$-name
of $\name \id_{\name{\bold d}}$.
\end{definition}

\begin{claim}
\label{c35}
For a frame $\bold d,\Vdash_{\bbQ_{\bar D}}
``\name \id_{\bold d}$ is an ideal on $\bbN$ containing the 
finite sets and $\bbN \notin \name \id_{\bold d}$"; moreover, for
every $A \in \cP(\bbN)$ from $\bold V$, we have $A = \emptyset \mod
F_{\bold d}$ iff $\Vdash_{\bbQ_{\bold d}} 
``A \in \name \id_{\bold d}"$.
\end{claim}

\begin{PROOF}{\ref{c35}}
It suffices to prove the following $\boxplus_1 - \boxplus_4$.
\mn
\begin{enumerate}
\item[$\boxplus_1$]  If $\Vdash_{\bbQ_{\bar D}}$ ``if $A_1 \subseteq
A_2$ and $A_2 \in \name \id_{\bold d}$ then $A_1 \in \name \id_{\bold d}"$.
\end{enumerate}
\mn
[Why?  If $(p,A,\bar S,\bar \zeta)$ is a local witness for 
$(\name A_1,\bold d)$ \then \, obviously it is a lcoal witness for 
$(\name A_2,\bold d)$.]
\mn
\begin{enumerate}
\item[$\boxplus_2$]  if $\Vdash_{\bbQ_{\bold d}}$ ``if
$A_1,A_2 \in \name \id_{\bold d}$ then $A_1 \cup 
A_2 \in \name \id_{\bold d}$".
\end{enumerate}
\mn
Why?  It suffices to prove:
 if $\Vdash_{\bbQ_{\bold d}} ``\name A_1
\cup \name A_2 = \name A \subseteq \bbN"$ and $\name A$ is 
$\bold d$-positive \then \,
$\name A_\ell$ is $\bold d$-positive for some $\ell \in \{1,2\}$.  Let
$(p,A,\bar S,\bar \zeta)$ be a local witness for $(\name A,\bold d)$
and we shall prove that there are $\ell \in \{1,2\}$ and a local witness for
$(\tr(p),\name A_\ell,\bold d)$; by the ``dense" in Definition
\ref{c31}(4) this suffices.

For any $n \in A$ and $\nu \in S_n$ such that $\zeta_n(\nu) = 0$ we
choose $(\ell_{n,\nu},\zeta_{n,\nu},S_{n,\nu})$ such that:
\mn
\begin{enumerate}
\item[$(*)_{2.1}$]  $(a) \quad \ell_{\nu,n} \in \{1,2\}$
\sn
\item[${{}}$]  $(b) \quad (S_{n,\nu},\zeta_{n,\nu}) \in \wfst(p^{[\nu]},\bar
D_{\bold d})$
\sn
\item[${{}}$]  $(c) \quad$ if $\zeta_{n,\nu}(\rho) = 0$ so $\rho \in
S_{n,\nu}$ \then \, there is $q \in \bbQ_{\bold d}$ such that

\hskip25pt $p \le q,\tr(q) = \rho$ and 
$q \Vdash ``n \in \name A_{\ell_{n,\nu}}"$; let $q_{n,\rho}$ be such $q$.
\end{enumerate}
\mn
[Why $(\rho_{n,\nu},\zeta_{n,\nu},S_{n,\nu})$ exists?  We shall use \ref{c30};
that is for $\ell \in \{1,2\}$ let $\cY_{n,\nu,\ell} = \{\rho:\nu
\trianglelefteq \rho \in p$ and there is $r \in \bbQ_{\bar D}$ such
that $\tr(r)=\rho$ and $p \le r$ and $r \Vdash ``n \in \name
A_\ell"\}$.

We apply for $\ell=1,2$ 
Claim \ref{c30} with $\bar D_{\bold d},\nu,\cY_{n,\nu,\ell}$
here standing for $\bar D,\eta^*,\cY$ there.  If for some $\ell \in
\{1,2\}$ clause (b) there holds as witness by the function $\zeta$,
easily the desired $(*)_{2.1}$ holds.  If for both $\ell = 1,2$ clause
(a) there holds then for $\ell = 1,2$ there is $q_\ell \in \bbQ_{\bold
d}$ such that $\tr(q_\ell) = \nu$ and $q_\ell \cap \cY_{n,\nu,\ell} =
\emptyset$.

Necessarily $q := q_1 \cap q_2 \cap p$ belongs to $\bbQ_{\bold d}$ and
has trunk $\nu$ and is disjoint to $\cY_{n,\nu,1} \cup
\cY_{n,\nu,2}$.  But $\bbQ_{\bold d} \models ``p^{[\nu]} \le q"$ and
$q^{[\nu]} \Vdash_{\bbQ_{\bold d}} ``n \in \name A = \name A_1 \cup
\name A_2"$, hence there are $\ell \in \{1,2\}$ and $r \in \bbQ_{\bold
d}$ such that $q \le r$ and $r \Vdash_{\bbQ_{\bold d}} ``n \in \name
A_\ell"$, but then $\tr(r) \in \cY_{n,\nu,\ell}$ and $\tr(r)  \in q_*
\subseteq q_\ell$, contradicting the choice of $q_\ell$.  So
$(*)_{2.1}$ holds indeed.]
\mn
\begin{enumerate}
\item[$(*)_{2.2}$]  \wilog \, $\zeta_{n,\nu}(\rho) = 0 \Rightarrow
\ell g(\rho) > n$.
\end{enumerate}
\mn
[Why?  Obvious.]
\mn
\begin{enumerate}
\item[$(*)_{2.3}$]  for $n \in A$ there are $\ell_n,S'_n,\zeta'_n$
such that
\sn
\begin{enumerate}
\item[$(a)$]  $(S'_n,\zeta'_n) \in \wfst(p,\bar D_{\bold d})$
\sn
\item[$(b)$]  $S'_n \subseteq S_n$ and $\max(S'_n) = S'_n \cap
\max(S_n)$
\sn
\item[$(c)$]  $\ell_n \in \{1,2\}$ and $\nu \in \max(S'_n) \Rightarrow
\ell_{n,\nu} = \ell_n$.
\end{enumerate}
\end{enumerate}
\mn
[Why?  Easy.]
\mn
\begin{enumerate}
\item[$(*)_{2.4}$]  for $n \in A$ letting $S''_n = \cup\{S_{n,\nu}:\nu
\in \max(S'_n)\} \cup S'_n$, for some $\zeta''_n$ and $\bar q_n$ we
have:
\sn
\begin{enumerate}
\item[$\bullet$]  $(S''_n,\zeta''_n) \in \wfst(p,\bar D)$
\sn
\item[$\bullet$]  $\{\rho:\zeta''_n(\rho)=0\} = \{\rho$: for some
$\nu$ we have $\nu \in S'_n,\zeta_n(\nu)=0,\rho \in S_{n,\nu}$ 

\hskip85pt and $\zeta_{n,\nu}(\rho)=0\}$ 
\sn
\item[$\bullet$]  $\bar q = \langle q_{n,\rho}:\zeta''_n(\rho) =
0\rangle$
\sn
\item[$\bullet$]  $\zeta''_n(\rho) = 0 \Rightarrow p \le q_{n,\rho}$
\sn
\item[$\bullet$]  $\tr(q_{n,\rho}) = \rho$
\sn
\item[$\bullet$]  $q_{n,\rho} \Vdash ``n \in \name A_{\ell_n}"$
\end{enumerate}
\end{enumerate}
\mn
[Why?  Think.]
\mn
\begin{enumerate}
\item[$(*)_{2.5}$]  there is $\ell \in \{1,2\}$ such that $A' := \{n
\in A:\ell_n=\ell\} \ne \emptyset \mod F_{\bold d}$.
\end{enumerate}
\mn
[Why?  Obvious as $A \in F^+_{\bold d}$.]

We now consider the quadruple $(p',A',\bar S'',\bar \zeta'')$ defined by:
\mn
\begin{enumerate}
\item[$\bullet$]  $p' = \{\varrho \in p$: if $\tr(p) \trianglelefteq
\rho \trianglelefteq \varrho,n \le \ell g(\rho)$ and $\rho \in \max(S''_n)$
then $\varrho \in q_{n,\rho}\}$ where $S''_n,q_{n,\rho}$ are from
$(*)_{2.4}$.
\end{enumerate}
\mn
[Why $p' \in \bbQ_{\bold d}$ with $\tr(p') = \tr(p)$?
  Recall$(*)_{2.2}$.]  

So together we have:
\mn
\begin{enumerate}
\item[$\bullet$]  $A'$ is from $(*)_{2.5}$, so $A' \in F^+_{\bold d}$
\sn
\item[$\bullet$]  $\bar S'' = \langle S''_n:n \in A'\rangle$ where
$S''_n$ is from $(*)_{2.4}$
\sn
\item[$\bullet$]  $\bar\zeta'' = \langle \zeta''_n:n \in A'\rangle$ where
$\zeta''_n$ is from $(*)_{2.4}$.
\end{enumerate}
\mn
Now check that $(p',A',\bar S'',\bar\zeta'')$ is a witness for
$(\tr(p),\name A_\ell,\bar D)$ hence $\boxplus_2$ holds as said in the
beginning of its proof.  
\mn
\begin{enumerate}
\item[$\boxplus_3$]   $\Vdash_{\bbQ_{\bold d}} ``\emptyset \in \name
\id_{\bold d}$; moreover if $A = \emptyset \mod F_{\bold d}$ is from
$\bold V$ then $A \in \name \id_{\bold d}"$.
\end{enumerate}
\mn
Why?  Because of clause (b) in Definition \ref{c31}(4).
\mn
\begin{enumerate}
\item[$\boxplus_4$]   $\Vdash_{\bbQ_{\bar D}[\bold d]} ``\bbN \notin 
\name \id_{\bold d}$, moreover if $B \in F^+_{\bold d}$ and $B \in \bold
V$ then $B \notin \name \id_{\bold d}"$.
\end{enumerate}
\mn
Why?  This means that $B$ is $\bold d$-positive which is obvious:
use the witness $(p,A,\bar S,\bar\zeta)$ where $p$ is any member of 
$\bbQ_{\bold d},A = B,S_n = \{\tr(p)\},\zeta_n(\tr(p)) = 0$.
\end{PROOF}

\begin{observation}
\label{c36}
Assume $\bold d_1,\bold d_2$ are frames and $\bar D_{\bold d_1} = \bar
D = \bar D_{\bold d_2}$ and $F_{\bold d_1} \subseteq F_{\bold d_2}$
\then \, $\Vdash_{\bbQ_{\bar D}} ``\name\id_{\bold d_1}
\subseteq \name \id_{\bold d_2}"$.
\end{observation}

\begin{PROOF}{\ref{c36}}
Should be clear.
\end{PROOF}

\begin{claim}
\label{c37}
We have $\Vdash_{\bbP_2} ``\name\id_{\name{\bold d}_1} \subseteq
\name\id_{\name{\bold d}_2}$ and $(\name\id_{\name{\bold
d}_1})^+[\bbP_1] \subseteq (\name\id_{\name{\bold d}_2})^+[\bbP_2]"$ \when \,:
\mn
\begin{enumerate}
\item[$(a)$]  $\bbP_1 \lessdot \bbP_2$
\sn
\item[$(b)$]  $\Vdash_{\bbP_\ell} ``\name{\bold d}_\ell$ is a frame"
for $\ell=1,2$
\sn
\item[$(c)$]  $\Vdash_{\bbP_2} ``\name D_{\bold d_1,\eta} \subseteq
\name D_{\bold d_2,\eta}"$ for $\eta \in {}^{\omega >}\omega$
\sn
\item[$(d)$]  if $A \in (\name D^+_{\name{\bold d}_1,\eta})^{\bold
V[\bbP_1]}$ then $A \in (\name D^+_{\name{\bold d}_2})^{\bold
V[\bbP_2]}$
\sn
\item[$(e)$]  $\Vdash_{\bbP_2} ``F_{\name{\bold d}_1} 
\subseteq F_{\name{\bold d}_2}"$
\sn
\item[$(f)$]  if $A \in (F^+_{\name{\bold d}_1})^{\bold V[\bbP_1]}$
then $A \in (F^+_{\name{\bold d}_2})^{\bold V[\bbP_2]}$.
\end{enumerate}
\end{claim}

\begin{PROOF}{\ref{c37}}
Should be clear by \ref{c39} below recalling \ref{cn.28}.
\end{PROOF}

\begin{claim}
\label{c39}
Let $\bold d$ be a frame and $\name A$ a $\bbQ_{\bar D_{\bold d}}$-name of a
subset of $\bbN$.  We have $\name A$ is $\bold d$-null \Iff \, for a
pre-dense set of $p \in \bbQ_{\bold d}$ we have $\tr(p)
\trianglelefteq \rho \in p \Rightarrow$ there is no local witness for
$(p^{[\rho]},\name A,\bold d)$ equivalently, for $(\rho,\name A,\bold d)$.
\end{claim}

\begin{PROOF}{\ref{c39}}
Straight.
\end{PROOF}

\begin{remark}
\label{c40}
The point of \ref{c39} is that the second condition is clearly
absolute in the relevant cases by \ref{cn.28}, i.e. in \ref{c37}.
\end{remark}

\begin{definition}
\label{c41}
1) fin$(I)$ is the set of finite functions from $I$ to
   $\cH(\aleph_0)$.

\noindent
2) Let $\bold K$ be the set of forcing notions $\bbQ$ such that some
pair $(I,f)$ witness it, i.e. $(I,f,\bbQ) \in \bold K^+$ which means:
\mn
\begin{enumerate}
\item[$(a)$]  $f$ is a function from $\bbQ$ to fin$(I)$
\sn
\item[$(b)$]  if $p_1,p_2 \in \bbQ$ and the functions $g(p_1),g(p_2)$
are compatible \then \, $p_1,p_2$ have a common upper bound $p$ with
$g(p) = g(p_1) \cup g(p_2)$.
\end{enumerate}
\mn
2) We define $\le_{\bold K} = \le^{\wk}_{\bold K}$ by:
$(I_1,f_1,\bbQ_1) \le^{\wk}_{\bold K} (I_2,f_2,\bbQ_2)$ means that:
\mn
\begin{enumerate}
\item[$(a)$]  $(I_\ell,f_\ell)$ witness $\bbQ_\ell \in \bold K$ for
$\ell=1,2$
\sn
\item[$(b)$]  $I_1 \subseteq I_2$
\sn
\item[$(c)$]  $f_1 \subseteq f_2$
\sn
\item[$(d)$]  $\bbQ_1 \subseteq_{\ic} \bbQ_2$.
\end{enumerate}
\mn
3) We define $\le^{\st}_{\bold K}$ similarly adding:
\mn
\begin{enumerate}
\item[$(d)^+$]  $\bbQ_1 \lessdot \bbQ_2$.
\end{enumerate}
\mn
4) If $\bold q \in \bold K^+$ let $\bold q = (I_{\bold q},f_{\bold
q},\bbQ_{\bold q})$.
\end{definition}

\begin{remark}
We can use much less in Definition \ref{c41}.
\end{remark}
\newpage

\section {Consistency of many gaps}

We prove the first result promised in the introduction.  Assume
$\lambda = \lambda^{< \lambda} > \aleph_1$ and we like to build a
c.c.c. forcing notion $\bbP$ of cardinality $\lambda$, such that
$\bold V^{\bbP}$ is as required: $\Sp_\chi$ includes $\Theta_1$ and is
disjoint to $\Theta_2$; really we force by $\bbP \times
\prod\limits_{\theta} \cT_\theta$, the $\cT_\theta$ quite complete and
translate $\bbP$-names of ultra systems of filters to ultra-filters.  
In order to have $\Theta_1 \subseteq
\Sp_\chi$, we shall represent $\bbP$ as an FS iteration $\langle
P_\alpha,\name{\bbQ}_\beta:\alpha \le \delta,\beta <
\delta\rangle,|\bbP_\alpha| \le \lambda$ and $\cT_\theta$ is,
e.g. ${}^{\theta >}2$ and for each $\theta \in \Theta_2$ we have a
$\bar D_\alpha = \langle \name D_{\alpha,s}:s \in \cT_\theta\rangle$ a
$\bbP_\alpha$-name of a ultra system of filters for unboundedly many
$\alpha < \delta$, increasing with $\alpha$; in the end we force by
$\bbP_\alpha \times \Pi\{\cT_\theta:\theta \in \Theta_1\}$.
Toward this for each $s \in \cT_\theta,\theta  \in \Theta_1$ we many
times force by $\bbQ^+_{D_{\alpha,s}}$ from \S1.  

But in order to have
$\Theta_2 \cap \Sp_\chi = \emptyset$, we intend to represent $\bbP$ as
the union of a $\lessdot$-increasing sequence $\langle
\bbP'_\varepsilon:\varepsilon < \lambda\rangle$ and for each $\theta
\in \Theta_2$ for stationarily many $\varepsilon <
\lambda,\cf(\varepsilon) = \theta$ and $\bbP'_{\varepsilon +1}$ is
essentially the ultrapower
$(\bbP'_\varepsilon)^\theta/E_\theta,E_\theta$ a $\theta$-complete
ultra-filter on $\theta$, so $\theta$ is a measurable cardinal.

To accomplish both we define a set $\bold Q$, each $\bold x \in \bold
Q$ consist of a FS-iteration of $\langle \bbP_\alpha,\bbQ_\beta:\alpha
\le \lambda^+,\beta < \lambda^+\rangle$ with $\langle \name{\bar
  D}_{s,\alpha}:s \in \cup\{\cT_\theta:\theta \in \Theta_1\}\rangle$
for many $\alpha < \lambda^+$, increasing with $\alpha$ and $\bbQ_\beta =
\bbQ_{D_{t(\beta),\alpha}}$. 

In the end for suitable $\bold x$, we shall use $\bbP_\delta$ for some
$\delta < \lambda^+$ of cofinality $\kappa << \lambda$, (e.g. $\kappa
= \aleph_1$). So why go so high as $\lambda^+$?  
It helps in the construction toward
the other aim; we shall construct $\langle \bold
x_\varepsilon:\varepsilon \le \lambda\rangle$ increasing in $\bold Q$
such that for each $\theta \in \Theta_2$ for $\varepsilon < \lambda$
of cofinality $\theta,\bold x_{\varepsilon +1}$ is essentially $(\bold
x_\varepsilon)^\theta/E_\theta$.  In particular, we have to prove
$\bold Q \ne \emptyset$, the existence of the ultrapower and the existence
of limit which happens to be a major proof here.  For this we have to
choose the right definition, in particular using
$\name{\text{id}}_{(D_{\alpha,s},D_{\beta,t})}$ from Definition \ref{c33}.

For this section we assume
\begin{hypothesis}
\label{cn.42}  
1) We now fix two cardinals $\kappa$ and $\lambda$ as well as two sets,
$\Theta_1$ and $\Theta_2$, of regular cardinals in the interval
$[\kappa,\lambda]$ and let $\Theta = \Theta_1 \cup \Theta_2$.  

Our assumptions are
\mn
\begin{enumerate}
\item[$(a)$]   $\kappa$ is regular and uncountable, $\lambda =
\lambda^{\aleph_0}$ and $\kappa < \lambda$
\sn
\item[$(b)$]   $\Theta_1$ and $\Theta_2$ are disjoint sets of regular
cardinals $< \lambda$ from the interval $[\kappa,\lambda)$ but $\kappa
\notin \Theta_2$
\sn
\item[$(c)$]  Each $\theta \in \Theta_1$ we have $\theta^{< \theta} =
\theta$
\sn
\item[$(d)$]   each $\theta \in \Theta_2$ carries a normal ultrafilter
$E_\theta$, hence $\Theta_2$ consists of measurable cardinals
\sn
\item[$(e)$]  for all $\theta \in \Theta_2$ the cardinal $\lambda$
satisfies $\cf(\lambda) > \theta$ and $\lambda =
\lambda^\theta/E_\theta$.
\end{enumerate}
\mn
2) Furthermore (and see \ref{cn.47}) below so it is not a burden
\mn
\begin{enumerate}
\item[$(f)$]   $\bar{\cT} = \langle {\cT}_\theta:\theta \in
\Theta_1\rangle,\cT_\theta$ is a tree of cardinality $\theta$
with $\theta$ levels, such
that above any element there are elements of any higher level (may add
``$\cT_\theta$ is $\aleph_2$-complete" and even ``$\cT_\theta$ is
$\theta$-complete", then clause (g) follows)
\sn
\item[$(g)$]   for every $\partial \in \Theta_1$, forcing by 
${\cT}_{\ge\partial} := \Pi\{{\cT}_\theta:\theta \in \Theta_1 \backslash
  \partial\}$, the product with Easton support, adds no 
sequence of ordinals of length $< \partial$ and, for
 simplicity, collapses no cardinal and changes no cofinality; if
$\kappa = \aleph_2 \in \Theta$ add ``$\cT_\kappa$ is
$\aleph_1$-complete"; let $\cT_* = \cT_{\ge \min(\Theta_1)}$
\sn
\item[$(h)$]  if $\partial \in \Theta_1$ then
$|\cT_{\ge \partial}|$ is $\Pi(\Theta_1 \backslash \partial)$ 
except when sup$(\Theta_1)$ is strongly  inaccessible 
and then the value is sup$(\Theta_1)$
\end{enumerate}
\end{hypothesis}

\begin{choice}
\label{c45}
1) Without loss of generality 
$\langle \cT_\theta:\theta \in \Theta_1 \rangle$ is a sequence of
 pairwise disjoint trees.

\noindent
2) Let $\cT$ be the disjoint sum of $\{\cT_\theta:\theta \in \Theta\}$, 
so it is a forest.

\noindent
3) Let $\bar t = \langle t_i:i \in S\rangle$ be a sequence of members of
$\cT$ where $S = \{\delta < \lambda^+:\cf(\delta) = \cf(\lambda)\}$ 
 such that if $t \in \cT$ then $\{\delta \in S:t_\delta = t\}$ 
is a stationary subset of $\lambda^+$; let $t(i) = t_i$.

\noindent
4) Furthermore choose
\mn
\begin{enumerate}
\item[$(\alpha)$]  $S_0 \subseteq \{\delta < \lambda^+:
\cf(\delta) = \aleph_0\}$ is stationary
\sn
\item[$(\beta)$]  $\bar\Upsilon = \langle
\Upsilon_{\delta,t,n}:\delta \in S_0,t \in \cT,n \in \bbN \rangle$; let
$\Upsilon(\delta,t,n) = \Upsilon_{\delta,t,n}$ 
\sn
\item[$(\gamma)$]  $\langle \Upsilon_{\delta,t,n}:n <
\omega\rangle$ is an increasing $\omega$-sequence of ordinals 
with limit $\delta$
\sn
\item[$(\delta)$]  $\Upsilon_{\delta,t,n} \in \{\alpha \in S:t_\alpha =t\}$
\sn
\item[$(\varepsilon)$]  $\bar\Upsilon$ guess clubs, i.e. if
$E$ is a club of $\lambda^+$ then the set
$\{\delta \in S_0:C^*_\delta :=
\{\Upsilon_{\delta,t,n}:t,n\} \subseteq E\}$ is stationary.
\end{enumerate}
\end{choice}

\begin{remark}
If $|\cT| <\lambda$ we can find such $\bar\Upsilon$, but in general it
is easy to force such $\bar\Upsilon$.
\end{remark}

\begin{claim}
\label{cn.47}
Assuming \ref{cn.42}(1) only, a 
sequence $\bar{\cT}$ as in \ref{cn.42}, clauses (f),(g),(h) (and
also $\bar t,s,S_\theta,\bar\Upsilon$ as in \ref{c45}) exists, 
provided that $\Theta_1 \subseteq
\{\theta:\theta = \theta^{< \theta} \ge \kappa\}$ and G.C.H. holds (or
just $\theta = \sup(\Theta_1 \cap \theta) \Rightarrow 2^\theta = \theta^+$).
\end{claim}

\begin{PROOF}{\ref{cn.47}}
Straight, e.g. $\cT_\theta = ({}^{\theta >}2,\triangleleft)$.
\end{PROOF}

\begin{definition}
\label{c49}  
Let $\bold Q$ be the set of objects $\bold x$ consisting of (below
$\alpha,\beta \le \lambda^+$):
\mn
\begin{enumerate}
\item[$(a)$]   $\bbP_\alpha \in \cH(\lambda^{++})$ and $I_{<
\alpha},f_\alpha \in \cH(\lambda^{++})$ 
witnessed $\bbP_\alpha \in \bold K$ for $\alpha \le 
\lambda^+$, all in $\cH(\lambda^+)$ if $\alpha <\lambda^+$
\sn
\item[$(b)$]   $I_\alpha \in \cH(\lambda^+)$ and
$\name{\bbQ}_\alpha,\name g_\alpha \in \cH(\lambda^+)$ are
$\bbP_\alpha$-names such that $\Vdash_{\bbP_\alpha} ``\name \bbQ_\alpha \in
\bold K$ as witnessed by $I_\alpha,\name g_\alpha"$ for $\alpha < \lambda^+$
\sn
\item[$(c)$]  $\bar{\bbQ} = \langle
\bbP_\alpha,\name{\bbQ}_\alpha:\alpha < \lambda^+ \rangle \in
\cH(\lambda^{++})$ is an FS iteration except that
\begin{enumerate}
\item[$(*)$]  $\bbP_\alpha = \{p:p$ a finite function with domain
$\subseteq \alpha$ such that if $\beta \in \dom(p)$ then
$\name g_\beta(p(\beta)) \in \fin(I_\beta)$ is an 
object (not just a $\bbP_\beta$-name)$\}$
\end{enumerate}
\item[$(d)$]   $I_{< \alpha} = \cup\{I_\beta:\beta < \alpha\}$ is
  disjoint to $I_\alpha$ and $\bbP =
\bbP_{\lambda^+} = \cup\{\bbP_\alpha:\alpha < \lambda^+\}$ and
$f_\alpha(\beta) = \cup\{\name g_\alpha(p(\beta)):\beta \in \dom(p)\}$
\sn
\item[$(e)$]   $E$ is a club of $\lambda^+$ and for $\alpha \in S \cap E$:
\begin{enumerate}
\item[$(\alpha)$]  $\name{\bar D}_\alpha = \langle
\name D_{\alpha,s}:s \in \cT\rangle$ is a $\bbP_\alpha$-name of an
ultra $\cT$-filter system (equivalently each $\name{\bar D}_{\alpha,\theta}
= \name{\bar D}_\alpha \rest \cT_\theta$ is a $\bbP_\alpha$-name of an ultra
$\cT_\theta$-filter system), and for simplicity 
$\fil(\name D_{\alpha,s}) = \name D_{\alpha,s}$
\sn
\item[$(\beta)$]  $\langle \name D_{\beta,s}:\beta \in S \cap E,\beta \le
\alpha\rangle$ is $\subseteq$-increasing continuous for each $s \in \cT$
\end{enumerate}
\item[$(f)$]   if $\alpha \in S \cap E$ \then \, 
$\bbQ_\alpha$ is $\bbQ_{\name D_{\alpha,t(\alpha)}}$
see Definition \ref{c20} and calling the generic 
$\name \eta_\alpha$, we have $I_\alpha = \{0\},\name g_\alpha(p) = \tr(p)$
\sn
\item[$(g)$]  $(\alpha) \quad$ if $\alpha \in S \cap E$ and $s,t \in
\cT$ \then \, $\Vdash_{\bbP_\alpha} \fil(\name D_{\alpha,s}) \subseteq
\fil(\name D_{\alpha,t})$ iff $s \le_{\cT} t$ actually
follows from $(e)(\alpha)$
\sn
\item[${{}}$]  $(\beta) \quad$ if $\alpha < \beta$ are from $S \cap E$
and $s \in \cT$ \then \, $\Vdash_{\bbP_\beta}$ ``if $\name A \in
\name\id(\name{\bold d}^\alpha_{t(\alpha),s})[\bbP_\alpha]$; 

\hskip25pt then $\name A = \emptyset \mod \name D_{\beta,s}"$ 
where $\name{\bold d}^\alpha_{t,s} = (\name D_{\alpha,t},\name D_{\alpha,s})$
\sn
\item[$(h)$]  $(\alpha) \quad$ if $\delta \in S_0 \cap E$ \then \,
$\bbQ_\delta$ is $\bbQ_{\fil(\emptyset)}$ with $\name \nu^*_\delta$
the generic
\sn
\item[${{}}$]  $(\beta) \quad$ if $\delta \in S_0 \cap E$ and $C^*_\delta
\subseteq E_{\bold x}$, see \ref{c45}$(4)(\varepsilon)$ then 

\hskip25pt $\name u_{\delta,t,n} \in \name D_{\gamma,t}$, see below, 
whenever $t \in \cT,n \in \bbN$ and 

\hskip25pt $\gamma \in S \cap E_{\bold x} \backslash (\delta +1)$
\sn
\item[${{}}$]  $(\gamma) \quad$ in clause $(\beta)$ we let $\name
u_{\delta,t,m} = \{\name\eta_{\Upsilon(\delta,t,n)}(k):n \in \bbN,n \ge m$

\hskip25pt  and $k \ge \name \nu^*_\delta(n)\}$.
\end{enumerate}
\end{definition}

\begin{discussion}
1) Later we shall use an increasing continuous sequence $\langle \bold
x_\varepsilon:\varepsilon \le \lambda\rangle$.  
Where and how will cofinality $\kappa$ reappear?  Well,
we shall use $\bbP_{\delta(*)}[\bold x_\lambda]$ for some $\delta(*)
\in E_{\bold x_\lambda}$ of cofinality $\kappa$.  So why not replace
$\lambda^+$ by $\kappa$ above?  We have a problem in proving the
existence of a (canonical) upper bound to $\langle \bold
x_\varepsilon:\varepsilon < \delta\rangle$, specifically in finding
the $\name D_{\beta_i}$ in the proof of Claim \ref{cn.56},
i.e. completing an appropriate $\cT$-filter system to an ultra one, e.g.
in Case 3 in the proof of \ref{cn.56}.
To help we carry a strong induction hypothesis, see clause
$(i)(\gamma)\bullet_2$ in $\boxdot$ there and then first find an
$\bbR_{\beta_j,\lambda^+}[\bbP_{\beta_j}\bar{\bold x}]$-name, then
reflect it to a $\beta_i$.

\noindent
2) Note that it helps to have not only $\name{\bbQ}_\alpha =
   \bbQ_{\bar D}$, but possibly some related forcing notions.  First
   in proving there is a limit, see \ref{cn.56}, in proving the ``reflection"
   discussed above lead us to use some unions.  Second, using
   ultrapower by $E_\theta$, see \ref{cn.77}, for limit $\delta$ of
   cofinality $\theta$, the ultrapower naturally leads us to use some
   iterations.

\noindent
3) We may in \ref{cn.42} demand $\kappa \notin \Theta_1$, equivalently
   $\kappa < \min(\Theta)$, but let $\cT_\kappa$ be a singleton
   $\{t_*\}$ and $\cT$ is $\cT_{\ge \min(\Theta_1)} \cup \cT_\kappa$.
   In this case in \ref{cn.91} we get $\Vdash_{\bbP \times \cT_*}
   ``\{\kappa\} \cup \Theta_1 \subseteq \Sp_\chi"$.
\end{discussion}

\begin{definition}
\label{c50}
1) For $\bold x \in \bold Q$, of course we let $\bar{\bbQ}^{\bold x} =
\bar{\bbQ}_{\bold x} = 
\bar{\bbQ}[\bold x] = \bar{\bbQ},\bbP^{\bold x}_\alpha =
\bbP_\alpha[\bold x] = \bbP_\alpha,\bbP_{\bold x} = \bbP^{\bold x} =
\bbP = \bbP^{\bold x}_{\lambda^+}$, etc.

\noindent
2) We define a two-place relation $\le_{\bold Q}$ on $\bold Q:\bold x
   \le_{\bold Q} \bold y$ iff:
\mn
\begin{enumerate}
\item[$(a)$]   $(I^{\bold x}_{< \alpha},f^{\bold x}_\alpha,\bbP^{\bold
x}_\alpha) \le^{\st}_{\bold K} (I^{\bold y}_{< \alpha},f^{\bold
y}_\alpha,\bbP^{\bold y}_\alpha)$ for $\alpha \le \lambda^+$, see
Definition \ref{c41}(3)
\sn
\item[$(b)$]  $\Vdash_{\bbP^{\bold y}_\alpha} ``(I^{\bold x}_\alpha,\name
g^{\bold x}_\alpha,\name{\bbQ}^{\bold x}_\alpha) \le^{\wk}_{\bold K}
(I^{\bold y}_\alpha,\name g^{\bold y}_\alpha,\name{\bbQ}^{\bold
y}_\alpha)"$ for $\alpha < \lambda^+$, see Definition \ref{c41}(2)
\sn
\item[$(c)$]  $E_{\bold y} \subseteq E_{\bold x}$
\sn
\item[$(d)$]  $\Vdash_{\bbP_\alpha[\bold y]} ``\name D^{\bold
x}_{\alpha,t(i)} \subseteq \name D^{\bold y}_{\alpha,t(i)}"$ for
$\alpha \in S \cap E_{\bold y}$ and $t \in \cT$
\sn
\item[$(e)$]  $\Vdash_{\bbP_\alpha[\bold y]}$ ``if $A \in ((\name
D^{\bold x}_{\alpha,t(i)})^+)^{\bold V[\bbP[\bold x]]}$ then $A \in
(\name D^{\bold y}_{\alpha,t(i)})^+"$, really follows by clause (d) and
\ref{c49}$(e)(\alpha)$, the ``ultra".
\end{enumerate}
\end{definition}

\begin{claim}
\label{cn.51}
$\bold Q$ is non-empty, in fact there is $\bold x \in \bold Q$ such
that $\bbP^{\bold x}_\alpha$ has cardinality $\lambda$ for $\alpha \in
[1,\lambda^+)$  and in $\bold V^{\bbP^{\bold x}_1}$ 
we have $2^{\aleph_0} = \lambda$.
\end{claim}

\begin{PROOF}{\ref{cn.51}}
For $i=0$, first letting $D'_{0,s} = \emptyset$ for $s \in \cT$,
clearly $\name D'_0 = \langle D'_{0,s}:s \in \cT\rangle$ is a 
$\cT$-filter system hence by \ref{cn.14}(2) we can choose $\bar D_0 =
\langle D_{0,s}:s \in \cT\rangle$, an ultra $\cT$-filter system (in $\bold V
= \bold V^{\bbP_0}$).  Second, we choose $\bbQ_i$ as adding $\lambda$ Cohen 
reals, say $\langle \name\eta^\ell_\alpha:\alpha < \lambda\rangle$ so
$I_i = \lambda,f_i$ is the identity, so
$f_i(p)(\alpha) = p(\alpha) \in {}^{\omega >}2$.  
Third, let $\langle (s_\alpha,t_\alpha):\alpha < \lambda\rangle$ be such that
$s_\alpha,t_\alpha \in \cT$ are $\le_{\cT}$-incomparable and any such
pair appears.

We define a $\bbP_1$-name $\name{\bar D}' = \langle \name D'_t:
t \in \cT \rangle$ by $\name D'_t =
\{\eta^{-1}_{1,\alpha}\{\ell\}:
s_\alpha \le_I t \wedge \ell = 0$ or
$t_\alpha \le_I t \wedge \ell=1\} \cup D_{0,t}$.  Clearly $\Vdash_{\bbP_1}
``\name{\bar D}'$ is an $\cT$-filter system", so by \ref{cn.14}(2)
there is $\name{\bar D}_1$ such that $\Vdash_{\bbP_1} 
``\name{\bar D}_1$ is an ultra $\cT$-filter satisfying 
$\name{\bar D}' \le \name{\bar D}_1$ hence $\bar D_0 \le \name{\bar D}_1$".

Now we shall choose $\bbP_\alpha,\name{\bar D}_\alpha$ by induction on
$\alpha \le \lambda^+$ also for $\alpha \in \lambda \backslash S$
 such that the relevant demands from Definition
\ref{c49} hold, in particular, $\langle \bbP_\beta,\name{\bbQ}_\gamma:
\beta \le \alpha,\gamma < \alpha\rangle$
is a FS-iteration but $\gamma \in \dom(p),p \in \bbP_\beta$ implies
that $\emptyset \in \bbP_\beta$ forces a value to $\tr(p(\gamma))$ and
also $\Vdash_{\bbP_\alpha} ``\name{\bar D}_\alpha$ is a
$\cT$-filter system such that $\name{\bar D}_\beta \le 
\name{\bar D}_\alpha$ for $\beta < \alpha$ and $\name{\bar D}_\alpha$
is ultra when $\alpha \notin S_0"$; recall that in Definition
\ref{c49} $\name{\bar D}_\alpha$ is defined only for $\alpha \in E
\cap S$, but no harm in defining $\name{\bar D}_\alpha$ in more cases.
For $\alpha = 0,1$ this was done above.

For $\alpha$ limit let $\bbP_\alpha = \cup\{\bbP_\beta:\beta <
\alpha\}$ and $\name{\bar D}'_\alpha = \langle \name D'_{\alpha,t}:t \in
\cT_*\rangle$ where 
$\name D'_{\alpha,t(i)} = \cup\{\name D_{\beta,t(i)}:\beta <
\alpha\}$.  It is easy to see that $\langle \bbP_\beta:\beta \le
\alpha\rangle$ is a $\lessdot$-increasing continuous sequence of
c.c.c. forcing notions and
$\Vdash_{\bbP_\alpha} ``\name{\bar D}'_\alpha$ is an $\cT$-filter
system".  If $\delta \in S_0$ let $\name{\bar D}_\alpha = \name{\bar
D}'_\alpha$, otherwise 
by \ref{cn.14}(2) we can find $\name{\bar D}_\alpha$ such
that $\Vdash_{\bbP_\alpha} ``\name{\bar D}_\alpha$ is an ultra
$\cT$-filter system and $\name{\bar D}'_\alpha \le \name{\bar D}_\alpha"$.

For $\alpha = \beta +1$ such that $\beta \notin S \cup S_0$ 
let $\name{\bbQ}_\beta$ be
trivial. Now let $\bbP_\alpha = \bbP_\beta * \name{\bbQ}_\beta$
and let $\name D'_{\alpha,t}$ be $\name D_{\beta,t}$.
Easily $\Vdash_{\bbP_\alpha} ``\langle \name D'_{\alpha,t}:t \in
\cT\rangle$ is an $\cT$-filter system" and choose 
$\name{\bar D}_\alpha$ as above, i.e. (a $\bbP_\alpha$-name of an) 
ultra $\cT$-filter system above $\name{\bar D}'_\alpha$.

Next, assume $\alpha = \beta +1,\beta \in S$; we let $\bbQ_\beta =
\bbQ_{\name D_{\beta,t(\beta)}}$ and $\bbP_\alpha = \bbP_\beta *
\name{\bbQ}_\beta$.
Now for $s \in \cT$, let $\name D'_{\alpha,s} = \name D_{\beta,s} \cup
\{\bbN \backslash A:A \in \name \id_{\name{\bold d}_{t(\beta),s}}\}$
where $t(\beta) = t_\beta$ is from \ref{c45}(3).
Note that $\Vdash_{\bbP_\beta} ``\fil(\name D_{\alpha,s}) \subseteq
\fil(\name D_{\alpha,t})$ iff $s \le_{\cT} t"$ by the choice of the
$\name D_{1,s}$'s and the $\name D_{\beta,s}$'s, so the definition of 
id$_{\name{\bold d}_{t(\beta),s}}$ depend on the 
truth value of $t(\beta) \le_I s$.

Now (pedantically working in $\bold V^{\bbP_\beta})$:
\mn
\begin{enumerate}
\item[$\bullet$]  $\name D'_{\alpha,s} \subseteq [\bbN]^{\aleph_0}$ by
its definition
\sn
\item[$\bullet$]  $ D_{\alpha,s} \subseteq D'_{\alpha,s}$, by \ref{c35}
\sn
\item[$\bullet$]  $\emptyset \notin \fil(D'_{\alpha,s})$ by \ref{c35}
\sn
\item[$\bullet$]  if $A \in (D^+_{\alpha,s})^{\bold V[\bbP_\beta]}$ then
$A \in ((D'_{\alpha,s})^+)^{\bold V[\bbP_{\beta +1}]}$ by \ref{c35}
\sn
\item[$\bullet$]   $s \le_I t \Rightarrow \name D'_{\alpha,s}
\subseteq \name D'_{\alpha,t}$ by \ref{c36} and the choice of the
$D'_{\alpha,t}$'s. 
\end{enumerate}
\mn
We continue as in the previous case.

Lastly, assume $\alpha = \beta +1,\beta \in S_0$ 
and we shall define for $\alpha$.
We let $\name{\bbQ}_\beta = \bbQ_{\fil(\emptyset)}$ in $\bold
V^{\bbP_\beta}$ and so $\name \nu^*_\beta$ is defined as the generic
and $\bbP_{\beta +1} = \bbP_\beta * \name{\bbQ}_\beta$.  Note that
$\name u_{\beta,t,n}$ is well defined,
(see clause (h) of Definition \ref{c49}).  By Claim
\ref{c52} below letting $\name D'_{\alpha,t} = \name D_{\beta,t}
\cup\{\name u_{\beta,s,n}:n \in \bbN$ and $s \in \cT$ satisfies 
$s \le_{\cT} t\}$ we have
$\name{\bar D}'_\alpha = \langle \name D'_{\alpha,t}:t \in \cT\rangle$
is a $\bbP_\beta$-name of a $\cT$-filter system above $\name{\bar
D}_\beta$ and let $\name{\bar D}_\alpha$ be (a $\bbP_\alpha$-name of)
an ultra $\cT$-filter system above $\name{\bar D}'_\alpha$.

Let $I_\alpha = \{\alpha\}$ for $\alpha < \lambda^+,
I_{< \alpha} = \alpha$ for $\alpha \le \lambda^+$ and if $\alpha \in S
\cup S_0$ then we let $\Vdash_{\bbP_\alpha}$
``if $p \in \name{\bbQ}_\alpha$ then 
$\name g_\alpha(p)$ is $\tr(p)$, the trunk" and if $\alpha \in
\lambda^+ \backslash (S \cup S_0)$ then $\name g_\alpha(p)=0$.

Naturally, we define $\bold x$ by: $\bbP^{\bold x}_\beta =
\bbP_\beta,\name{\bbQ}^{\bold x}_\alpha =
\name{\bbQ}_\alpha,E_{\bold x} = \lambda,I^{\bold x}_\alpha =
I_\alpha,I^{\bold x}_{< \beta} = I_{< \beta},\name g^{\bold x}_\alpha
= \name g_\alpha$ for $\alpha < \lambda^+,\beta \le \lambda^+$ (and so
$f^{\bold x}_\alpha$ is defined), 
$\name{\bar D}^{\bold x}_\gamma = \name D_\gamma$ for $\gamma \in S,
\beta \le \lambda^+,\alpha < \lambda$.  It is easily to check that $\bold
x \in \bold Q$ is as required.
\end{PROOF}

\begin{claim}
\label{c52}
If (A) then (B) where
\mn
\begin{enumerate}
\item[$(A)$]  $(a) \quad \delta \in S_0$
\sn
\item[${{}}$]  $(b) \quad \bbP_\alpha(\alpha \le
\delta),\name{\bbQ}_\alpha(\alpha \le \delta),E \subseteq \delta$,
etc., are as in Definition \ref{c49} except

\hskip25pt  that all is up to $\delta$
\sn
\item[${{}}$]  $(c) \quad \name{\bbQ}_\delta,\name\nu^*_\delta,\name
u_{\delta,t}$ are as in clause (h) of Definition \ref{c49}
\sn
\item[${{}}$]  $(d) \quad \name D'_{\delta,t} := \cup\{\name
D_{\alpha,t}:\alpha \in S \cap E\} \cup \{\name u_{\delta,t,n}:n \in
\bbN$ and $s \in \cT$ satisfies 

\hskip25pt  $s \le_{\cT} t\}$ so a 
$\bbP_\delta * \name{\bbQ}_{\fil(\emptyset)}$-name
\sn
\item[$(B)$]  $(a) \quad \Vdash_{\bbP_\delta *
\name{\bbQ}_{\fil(\emptyset)}} ``\langle \name D'_{\delta,t}:t \in
\cT\rangle$ is a $\cT$-filter system"
\sn
\item[${{}}$]  $(b) \quad \Vdash_{\bbP_\delta *
\name{\bbQ}_{\fil(\emptyset)}} ``\fil(\name D'_{\delta,t}) =
\fil(\{\name u_{\delta,s,n}:s \le_{\cT} t$ and $n \in \bbN\})"$
\sn
\item[${{}}$]  $(c) \quad  \Vdash_{\bbP_\delta *
\name{\bbQ}_{\fil(\emptyset)}}$ ``if $t \in \cT$ and
$A \in \cup\{D_{\alpha,t}:\alpha \in \delta \cap S\}$ then 
$\name u_{\delta,t,n} \subseteq^* A$ 

\hskip25pt  for every large enough $n"$.
\end{enumerate}
\end{claim}

\begin{PROOF}{\ref{c52}}
Straight; the point is $\Vdash_{\bbP_\delta *
\name{\bbQ}_{\fil(\emptyset)}} ``\emptyset \in \fil(\name
D'_{\delta,t})"$ for $t \in \cT$, which holds as
\mn
\begin{enumerate}
\item[$(*)_1$]   if $A \in \name D_{\Upsilon(\delta,t,n)}$
\then \, for every large enough $k,\name\eta_{\Upsilon(\delta,t,n)}(k) \in A$
\sn
\item[$(*)_2$]   if $A \in \name D^+_{\Upsilon(\delta,t,n)}$
in $\bold V^{\bbP_\delta}$ \then \, for infinitely many 
$k,\name\eta_{\Upsilon(\delta,t,n)}(k) \in A$
\sn
\item[$(*)_3$]  $\name\nu^*_\delta$ is a dominating real.
\end{enumerate}
\end{PROOF}

\begin{observation}
\label{cn.53} 
1) $\le_{\bold Q}$ partially orders $\bold Q$.

\noindent
2) $\bbP^{\bold x}_\alpha$ satisfies the c.c.c. and even is locally
   $\aleph_1$-centered\footnote{meaning that any $\aleph_1$ elements can be
   divided to $\aleph_0$ sets such that any finitely many members of
   one sets has a common upper bound} when $\bold x \in \bold Q$ and
   $\alpha \le \lambda^+$.
\end{observation}

\begin{PROOF}{\ref{cn.53}}
Easy.
\end{PROOF}

\begin{claim}
\label{cn.56}
\underline{The upper bound existence claim}  

If $\langle \bold x_\varepsilon:\varepsilon < \delta\rangle$ is
$\le_{\bold Q}$-increasing and $\delta$ is a limit ordinal $<
\lambda^+$ \then \, there is $\bold x_\delta$ which is a canonical
limit of $\langle \bold x_\varepsilon:\varepsilon < \delta\rangle$, see below.
\end{claim}

\begin{definition}
\label{cn.57}
We say $\bold x = \bold x_\delta$ is a canonical limit of $\bar{\bold
x} = \langle \bold x_\varepsilon:\varepsilon < \delta\rangle$ \when \,
$\bar{\bold x}$ is $\le_{\bold Q}$-increasing, $\delta$ is a limit
ordinal $< \lambda^+$ and (for every $\alpha < \lambda^+$):
\mn
\begin{enumerate}
\item[$(a)$]   $\bold x_\delta \in \bold Q$
\sn
\item[$(b)$]   $\bold x_\varepsilon \le_{\bold Q} \bold x_\delta$ for
$\varepsilon < \delta$ and $E_{\bold x_\delta} \subseteq
\cap\{E_{\bold x_\varepsilon}:\varepsilon < \delta\}$ 
\sn
\item[$(c)$]   $I_\alpha[\bold x_\delta] = \cup\{I_\alpha[\bold
x_\varepsilon]:\varepsilon < \delta\}$
\sn
\item[$(d)$]   if $\delta$ has uncountable cofinality \then \,
\begin{enumerate}
\item[$(\alpha)$]   $\bbP_\alpha^{\bold x_\delta} = 
\cup\{\bbP^{\bold x_\varepsilon}_\alpha:\varepsilon < \delta\}$ 
\sn
\item[$(\beta)$]    $\Vdash_{\bbP^{\bold x_\delta}_\alpha} 
``\name D^{\bold x_\delta}_{\alpha,t} =
\cup\{\name D^{\bold x_\varepsilon}_{\alpha,t}:\alpha < \delta\}"$ for
 $t \in \cT$ if $\alpha \in E_{\bold x_\delta} \cap S$
\sn
\item[$(\gamma)$]  $\name{\bbQ}^{\bold x_\delta}_\alpha =
\cup\{\name{\bbQ}^{\bold x_\varepsilon}_\alpha:\alpha < \delta\}$
\sn
\item[$(\delta)$]  $\name g^{\bold x_\varepsilon}_\alpha = \cup\{\name
g^{\bold x_\varepsilon}_\alpha:\varepsilon < \delta\}$.
\end{enumerate}
\item[$(e)$]   if $\delta$ has cofinality $\aleph_0$, then
\begin{enumerate}
\item[$(\alpha)$]  if $\alpha \in \lambda^+ \backslash 
(S \cap E_{\bold x_\delta})$ or $\alpha \in S_0 \wedge C^*_\alpha
\nsubseteq E_{\bold x_\delta}$ then $\Vdash_{\bbP_\alpha[\bold x_\delta]} 
``\name{\bbQ}_\alpha[\bold x_\delta] =
\cup\{\name{\bbQ}_\alpha[\bold x_\varepsilon]:\varepsilon < \delta\}$
and similarly $\name g_\alpha[\bold x_\delta] = \cup\{\name
g_\alpha[\bold x_\varepsilon]:\varepsilon < \delta\}$
\sn
\item[$(\beta)$]  if $\alpha \in S \cap E_{\bold x_\delta}$ then 
$\Vdash_{\bbP_\alpha[\bold x_\delta]} ``\name
D_{\alpha,t}[\bold x_\delta] \supseteq \cup\{\name D_{\alpha,t}[\bold
x_\varepsilon]:\varepsilon < \delta\}"$
\end{enumerate}
\item[$(f)$]  in fact $|\bbP^{\bold x_\delta}_\alpha| \le 
(\Sigma\{|\bbP^{\bold x_\varepsilon}_\alpha|:
\varepsilon < \delta\})^{\aleph_0}$.
\end{enumerate}
\end{definition}

\begin{PROOF}{\ref{cn.56}}
Let
\mn
\begin{enumerate}
\item[$\boxplus_0$]  $(a) \quad I_\alpha = 
\cup\{I_\alpha[\bold x_\varepsilon]:\varepsilon <
\delta\}$ for $\alpha < \lambda^+$
\sn
\item[${{}}$]  $(b) \quad I_{< \alpha} =
\cup\{I_\beta:\beta < \alpha\}$ for $\alpha \le \lambda^+$
\sn
\item[${{}}$]  $(c) \quad E := \cap\{E[\bold x_\varepsilon]:
\varepsilon < \delta\}$.
\end{enumerate}
\mn
So $E \subseteq \cap\{E[\bold x_\varepsilon]:\varepsilon <
\delta\}$ and 
clearly $E$ is a club of $\lambda^+$ (but in general this will not
be $E[\bold x_\delta]$).  If $\beta \le \gamma \le \lambda^+$ and $\bbQ$
satisfies $\varepsilon < \delta  \Rightarrow \bbP_\beta[\bold
x_\varepsilon] \lessdot \bbQ$ and for transparency $q \in \bbQ
\Rightarrow \emptyset \le_{\bbQ} q$ then $\bbR =
\bbR_{\beta,\gamma}[\bbQ,\bar{\bold x}]$ is defined as follows:
\mn
\begin{enumerate}
\item[$\boxplus_1$]  $(a) \quad p \in \bbR$ iff $p = (p_1,p_2)$
and some pair $(\varepsilon,p_0)$ witenss it which means

\hskip25pt $\varepsilon < \delta$ and $p_0 \in
\bbP_\beta[\bold x_\varepsilon],p_1 \in \bbP_\gamma[\bold
x_\varepsilon],p_2 \in \bbQ$ and one of the

\hskip25pt  following occurs
\begin{enumerate}
\item[${{}}$]  $(\alpha) \quad p_1 = \emptyset$ or $p_2 = \emptyset$
recalling clause (c) of \ref{c49}
\sn
\item[${{}}$]  $(\beta) \quad p_0 \Vdash_{\bbP_\beta[\bold
x_\varepsilon]} ``p_1 \in \bbP_\gamma[\bold x_\varepsilon]/
\bbP_\beta[\bold x_\varepsilon]$ and $p_2 \in 
\bbQ/\bbP_\beta[\bold x_\varepsilon]"$
\end{enumerate}
\item[${{}}$]  $(b) \quad$ for $p \in \bbR$ let $\varepsilon(p)$ be
the minimal $\varepsilon < \delta$ such that $(\varepsilon,p_0)$

\hskip25pt witness $p \in \bbR$ for some $p_0$
\sn
\item[${{}}$]  $(c) \quad \bbR \models ``p \le q"$ iff letting
$\varepsilon = \max\{\varepsilon(p),\varepsilon(q)\}$ we have

\hskip25pt $\bbP_\gamma[\bold x_\varepsilon] \models ``p_1 \le q_1"$ and
$\bbQ \models ``p_2 \le q_2"$.
\end{enumerate}
\mn
We note that:
\mn
\begin{enumerate}
\item[$\boxplus_2$]  $(a) \quad \bbR_{\beta,\gamma}[\bbQ,\bar{\bold
    x}]$ is a partial order
\sn
\item[${{}}$]  $(b) \quad$ above $\bbR'_{\beta,\gamma}[\bbQ,\bar{\bold x}]$
is a dense subset of $\bbR_{\beta,\gamma}[\bbQ,\bar{\bold x}]$ where
$\bbR'_{\beta,\gamma}[\bbQ,\bar{\bold x}]$ is defined like
$\bbR_{\beta,\gamma}[\bbQ,\bar{\bold x}]$ when in $\boxplus_1(a)$ we
omit subclause $(\alpha)$.
\end{enumerate}
\mn
[Why?  Clause (a) by $\boxplus_3$ below and clause (b) is easy.]

So below we may ignore the difference between
$\bbR_{\beta,\gamma}[\bbQ,\bar{\bold x}]$ and
$\bbR'_{\beta,\gamma}[\bbQ,\bar{\bold x}]$
\mn
\begin{enumerate}
\item[$\boxplus_3$]  for $(\beta,\gamma,\bbQ)$ as above; if
$(\varepsilon,p_0)$ is a witness for $p=(p_1,p_2) \in
\bbR_{\beta,\gamma}[\bbQ,\bar{\bold x}]$ and $\zeta \in
(\varepsilon,\delta)$ \then \, for some $q_0 \in \bbP_\beta[\bold
x_\zeta]$ the pair $(\zeta,q_0)$ is a witness for $(p_1,p_2) \in
\bbR_{\beta,\gamma}[\bbQ,\bar{\bold x}]$.
\end{enumerate}
\mn
[Why?  As we can increase $p_0$ in $\bbP_\beta[\bold x_\varepsilon]$,
\wilog \, $(p_1 \rest \beta) \le p_0$, where on $\rest$ recall Definition
\ref{c49}, clause (c).  As $(\varepsilon,p_0)$ is a witness for
$(p_1,p_2) \in \bbR_{\beta,\gamma}[\bbQ,\bar{\bold x}]$ necessarily
$p_0,p_2$ are compatible in $\bbQ$ hence they have a common upper
bound $q_2 \in \bbQ$.  As $\bbP_\beta[\bold x_\zeta] \lessdot \bbQ$, 
there is $q_0 \in \bbP_\beta[\bold x_\zeta]$ such that $q_0 \le q
\in \bbP_\beta[\bold x_\zeta] \Rightarrow q,q_2$ are compatible in
$\bbQ$.  As we can increase $q_0$ in $\bbP_\beta[\bold x_\zeta]$ and $p_0
\le q_2$ \wilog \, $p_0 \le q_0$ but $(p_1 \rest \beta) \le p_0$ hence
$(p_1 \rest \beta) \le q_0$.  As $\bold x_\varepsilon \le \bold
x_\zeta$ and $\langle \bbP_\alpha[\bold
x_\zeta],\name{\bbQ}_\alpha[\bold x_\zeta]:\alpha < \lambda^+\rangle$
is FS iteration and $p_1 \in \bbP_\gamma[\bold x_\varepsilon] \lessdot
\bbP_\gamma[\bold x_\zeta]$, clearly $q_0 \le q \in \bbP_\beta[\bold
x_\zeta] \Rightarrow q,p_1$ are compatible.  So clearly $(\zeta,q_0)$
is a witness for $p \in \bbR_{\beta,\gamma}[\bbQ,\bar{\bold x}]$ as
required in $\boxplus_3$.]
\mn
\begin{enumerate}
\item[$\boxplus_4$]  if $\beta,\gamma,\bbQ$ are as above and
$\gamma \le \gamma(1) \le \lambda^+$ \then \,
$\bbR_{\beta,\gamma}[\bbQ,\bar{\bold x}] \lessdot
\bbR_{\beta,\gamma(1)}[\bbQ,\bar{\bold x}]$.
\end{enumerate}
\mn
[Why?  We check the conditions from Definition \ref{z4}(3), the second
alternative.  First, if $p=(p_1,p_2) \in
\bbR_{\beta,\gamma}[\bbQ,\bar{\bold x}]$ we shall prove $p \in
\bbR_{\beta,\gamma(1)}[\bbQ,\bar{\bold x}]$; as $p \in
\bbR_{\beta,\gamma}[\bbQ,\bar{\bold x}]$, some $(\varepsilon,p_0)$
witness it, easily it witnesses $p \in
\bbR_{\beta,\gamma(1)}[\bbQ,\bar{\bold x}]$ as $\bbP_\gamma[\bold
x_\varepsilon] \subseteq \bbP_{\gamma(1)}[\bold x_\varepsilon]$.

Second, assume $\bbR_{\beta,\gamma}[\bbQ,\bar{\bold x}] \models 
``p \le q"$ and we should prove
$\bbR_{\beta,\gamma(1)}[\bbQ,\bar{\bold x}] \models ``p \le q"$, this
is obvious by the definition of the orders for those forcing notions.
Together $\bbR_{\beta,\gamma}[\bbQ,\bar{\bold x}] \subseteq
\bbR_{\beta,\gamma(1)}[\bbQ,\bar{\bold x}]$.

Third, we should prove $\bbR_{\beta,\gamma}[\bbQ,\bar{\bold x}]
\subseteq_{\ic} \bbR_{\beta,\gamma(1)}[\bbQ,\bar{\bold x}]$ so assume
$p,q \in \bbR_{\beta,\gamma}[\bbQ,\bar{\bold x}]$ has a common upper
bound $r=(r_1,r_2)$ in $\bbR_{\beta,\gamma(1)}[\bbQ,\bar{\bold x}]$.
Now easily $(r_1 \rest \gamma,r_2)$ is a common upper bound of $p,q$
in $\bbR_{\beta,\gamma}[\bbQ,\bar{\bold x}]$ as required.

Fourth, for $p \in \bbR_{\beta,\gamma(1)}[\bbQ,\bar{\bold x}]$ we
should find $q \in \bbR_{\beta,\gamma}[\bbQ,\bar{\bold x}]$ such that
if $\bbR_{\beta,\gamma}[\bbQ,\bar{\bold x}] \models ``q \le q^*"$ then
$q^*,p$ are compatible in $\bbR_{\beta,\gamma(1)}[\bbQ,\bar{\bold x}]$.

Now let $p = (p_1,p_2) \in \bbR_{\beta,\gamma(1)}[\bbQ,\bar{\bold x}]$ and let
$(\varepsilon,p_0)$ witness it; \wilog \, $\bbP_\beta[\bold
x_\varepsilon] \models ``(p_1 \rest \beta) \le p_0"$.

Let $q_1 = p_1 \rest \gamma \in \bbP_\gamma[\bold x_\varepsilon]$, now $q :=
(q_1,p_2)$ satisfies
\mn
\begin{enumerate}
\item[$\bullet$]  $q \in \bbR_{\beta,\gamma}[\bbQ,\bar{\bold x}]$.
\end{enumerate}
\mn
Why?  The pair $(\varepsilon,p_0)$ witness it because if $p_0 \le q'
\in \bbP_\beta[\bold x_\varepsilon]$ \then \, first $p_1,q'$ has a
common upper bound $r \in \bbP_{\gamma(1)[\bold x_\varepsilon]}$ hence
$r \rest \gamma \in \bbP_\gamma[\bold x_\varepsilon]$ is a common
upper bound of $q',q_1$; second $q',p_2$ has a common upper bound in
$\bbQ$ as $(\varepsilon,p_0)$ witness $(p_1,p_2)$.  So indeed
$(\varepsilon,p_0)$ witness $q = (q_1,p_2) \in
\bbR_{\beta,\gamma}[\bbQ,\bar{\bold x}]$.
\mn
\begin{enumerate}
\item[$\bullet$]  If $q \le q^* \in
\bbR_{\beta,\gamma}[\bbQ,\bar{\bold x}]$ \then \, $q^*,p$ are
compatible in $\bbP_{\gamma(1)}[\bold x_\varepsilon]$.
\end{enumerate}
\mn
Why?  Let $q^* = (q^*_1,q^*_2)$ and let $r_1 = (p_1 \rest
[\gamma,\gamma(1))) \cup q^*_1$, easily $(r_1,q^*_2) \in 
\bbR_{\beta,\gamma(1)}[\bbQ,\bar{\bold x}]$ is a common upper bound of
$q^*,p$.

This finishes checking the last demand for
``$\bbR_{\beta,\gamma}[\bbQ,\bar{\bold x}] \lessdot
\bbR_{\beta,\gamma(1)}[\bbQ,\bar{\bold x}]$ so $\boxplus_4$ holds.]
\mn
\begin{enumerate}
\item[$\boxplus_5$]  if $\bbQ$ satisfies the c.c.c. then 
$\bbR_{\beta,\gamma}[\bbQ,\bar{\bold x}]$ satisfies the c.c.c..
\end{enumerate}
\mn
[Why?  Let $p_i = (p_{1,i},p_{2,i}) \in \bbR_{\beta,\gamma}
[\bbQ,\bar{\bold x}]$ for $i < \aleph_1$.  Let
$(\varepsilon_i,p_{0,i})$ be a witness for $(p_{1,i},p_{2,i})$.  As
before let $q_i \in \bbQ$ be such that $p_0,p_{1,i} \rest \beta,p_{2,i}$
are below it.

We can find an uncountable $S$ such that $\langle f_\gamma[\bold
 x_{\varepsilon_i}](p_{1,i}):i \in S \rangle$ are pairwise compatible
 functions and $\langle \varepsilon_i:i \in S\rangle$ is
 non-decreasing.  As $\bbQ$ satisfies the c.c.c., for some $i<j$ from
 $S$ there is a common upper bound $q \in \bbQ$ of $q_i,q_j$; let
 $\{\beta_\ell:\ell < n\}$ list in increasing order $\{\beta\} \cup
 \dom(p_{1,i}) \cup \dom(p_{1,j}) \backslash \beta$ and let $\beta_n =
 \gamma$.

By induction on $\ell \le n$ we 
choose $r_\ell \in \bbP_{\beta_\ell}[\bold x_{\varepsilon_j}]$ such that:
\mn
\begin{enumerate}
\item[$\bullet$]  if $\ell=0$ so $\beta_\ell = \beta$ then $r_0 \le
r \in \bbP_\beta[\bold x_{\varepsilon_j}] \Rightarrow r,q$ are 
compatible in $\bbQ$
\sn
\item[$\bullet$]  if $\ell=m+1$ then $r_m \le r_\ell$
\sn
\item[$\bullet$]   $\bbP_{\beta_\ell}[\bold x_{\varepsilon_j}] \models
``(p_{1,i} \rest \beta_\ell) \le r_\ell$ and $(p_{1,j} \rest \beta_\ell)
\le r_\ell"$.
\end{enumerate}
\mn
For $\ell=0$ use $q \in \bbQ$ and
$\bbP_\beta[\bold x_{\varepsilon_j}] \lessdot \bbQ$.
For $\ell=m+1$, we shall choose $r_\ell \in
\bbP_{\beta_m+1}[\bold x_{\varepsilon_i}]$ as follows: if $\beta_\ell \notin
\dom(p_{1,i})$ then $r_\ell = r_m
\cup\{(\beta_\ell,p_{1,j}(\beta_\ell))\}$; if 
$\beta_\ell \notin \dom(p_{1,j})$ similarly;
otherwise, i.e. if $\beta_\ell \in \dom(p_{1,i}) \cap \dom(p_{1,j})$ 
use the demands on $\name g_{\beta_\ell}$ recalling $(*)$ of clause
(c) and end of clause (d) of Definition \ref{c49}.

Having carried the induction, $(r_m,q)$ is well defined.
Now let $r_* \in 
\bbP_\beta[\bold x_{\varepsilon_j}]$ be above $r_0$ such that $r_* \le
r \in \bbP_\beta[\bold x_{\varepsilon_j}] \Rightarrow r_m,r$ are
compatible.  Also $r_* \le r \in \bbP_\beta[\bold x_{\varepsilon_j}]
\Rightarrow r_0 \le r \in \bbP_\beta[\bold x_{\varepsilon_j}]
\Rightarrow r,q$ are compatible in $\bbQ$.  So $(\varepsilon_j,r_*)$
witness $(r_m,q) \in \bbR_{\beta,\gamma}[\bbQ,\bar{\bold x}]$ and
easily $(r_m,q)$ is above $p_i = (p_{1,i},p_{2,i})$ and above $p_j =
(p_{1,j},p_{2,j})$, so $\boxplus_5$ holds indeed.]
\mn
\begin{enumerate}
\item[$\boxplus_6$]  for $\beta,\gamma,\bbQ$ as above, $\bbQ \lessdot
\bbR_{\beta,\gamma}[\bbQ,\bar{\bold x}]$ when we identify $p_2 \in
\bbQ$ with $(\emptyset,p_2)$.
\end{enumerate}
\mn
[Why?  Again, first $p \in \bbQ \Rightarrow p \in 
\bbR_{\beta,\gamma}[\bbQ,\bar{\bold x}]$ by the identification, and
for $p,q \in \bbQ$ we have $\bbQ \models ``p \le q" \Leftrightarrow
\bbR_{\beta,\gamma}[\bbQ,\bar{\bold x}] \models ``p \le q"$ by the
definition of the order of $\bbR_{\beta,\gamma}[\bbQ,\bar{\bold x}]$.
So $\bbQ \subseteq \bbR_{\beta,\gamma}[\bbQ,\bar{\bold x}]$ holds,
moreover $\bbQ \subseteq_{\ic} \bbR_{\beta,\gamma}[\bbQ,\bar{\bold
x}]$ by the definition of the order.

Lastly, let $q \in \bbR_{\beta,\gamma}[\bbQ,\bar{\bold x}]$, so by
$\boxplus_2$ \wilog
\, $q = (q_1,q_2) \in \bbR'_{\beta,\gamma}[\bbQ,\bar{\bold
x}]$ and we shall find $p \in \bbQ$ such that $p \le p' \in \bbQ
\Rightarrow p',(q_1,q_2)$ are compatible.

Let $p=q_2$, i.e. $(\emptyset,q_2)$, and the rest should be clear.]
\mn
\begin{enumerate}
\item[$\boxplus_7$]  for $\beta,\gamma,\bbQ$ as above we have
$\bbP_\gamma[\bold x_\varepsilon] \lessdot \bbR_{\beta,\gamma}
[\bbQ,\bar{\bold x}]$ when we identify $p_1 \in \bbP_\gamma[\bold
x_\varepsilon]$ with $(p_1,\emptyset)$.
\end{enumerate}
\mn
[Why?  Similarly.]
\bigskip

\centerline {$* \qquad * \qquad *$}
\bigskip

Now by induction on $i \le \lambda^+$ we choose $\beta_i$ and
$\bbP_\alpha,f_\alpha$ 
(when $\alpha \le \beta_i$ and $j < i \Rightarrow \beta_j < \alpha$),
$\name{\bbQ}_\alpha,\name g_\alpha$ 
(when $\alpha < \beta_i$ and $j<i \Rightarrow
\beta_j \le \alpha)$ and\footnote{so we define some $\name{\bar
D}_{\beta_i}$ not used in $\bold x_\delta$.}
 also $\name{\bar D}_{\beta_i}$ 
(when $\beta_i \in S$) such that
\mn
\begin{enumerate}
\item[$\boxdot$]   the relevant parts of clauses (a)-(e) of Definition
\ref{cn.57} and of the definition of $\bold x_\delta \in \bold Q$ holds, in
particular (all when defined):
\begin{enumerate}
\item[$(a)$]  $\bbP_\alpha \in \cH(\lambda^+)$ is a c.c.c. forcing notion
\sn
\item[$(b)$]  $(\alpha) \quad \bbP^{{\bold x}_\varepsilon}_\alpha
\lessdot \bbP_\alpha$ and $\bbP_{\bold x_\varepsilon} \cap \bbP_\alpha
= \bbP^{{\bold x}_\varepsilon}_\alpha$ for $\varepsilon < \delta$
\sn
\item[${{}}$]  $(\beta) \quad (I_{< \alpha}[\bold
x_\varepsilon],f_\alpha[\bold x_\varepsilon],\bbP_\alpha[\bold
x_\varepsilon]) \le^{\st}_{\bold K} (I_{< \alpha},f_\alpha,\bbP_\alpha)$
\sn
\item[$(c)$]  $\name{\bar D}_{\beta_i}$
is a $\bbP_{\beta_i}$-name of an $I$-filter system; ultra when
$\beta_i \in S$; see $(i)(\gamma)\bullet_1$
\sn
\item[$(d)$]   if $\beta_i \in S,\varepsilon < \delta$ and $t \in \cT$ 
then $\Vdash_{\bbP_{\beta_i}} ``\name D^{\bold x_\varepsilon}_{\beta_i,t}
\subseteq \name D_{\beta_i,t}"$ 
\sn
\item[$(e)$]  $\langle \bbP_\alpha:\alpha \le \beta_i\rangle$ 
is $\lessdot$-increasing continuous
\sn
\item[$(f)$]   if $\beta = \alpha +1$ then $\bbP_\gamma 
= \bbP_\alpha * \name{\bbQ}_\alpha$, in fact, $\langle
\bbP_\beta,\name{\bbQ}_\alpha:\beta \le \beta_i,\alpha <
\beta_i\rangle$ is as in clause (c) of Definition \ref{c49}
\sn
\item[$(g)$]   if $\neg(\exists j)(\alpha = \beta_j \in S)$ 
then $\Vdash_{\bbP_\alpha} ``\name{\bbQ}_\alpha =
\cup\{\name{\bbQ}_\alpha[\bold x_\varepsilon]:\varepsilon <
\delta\},\name g_\alpha = \cup\{\name g_\alpha[\bold
x_\varepsilon]:\varepsilon < \delta\}"$; note that
$\name{\bbQ}_\alpha[\bold x_\varepsilon],\name g_\alpha[\bold
x_\varepsilon]$ are $\bbP^{\bold x_\varepsilon}_\alpha$-names hence
$\bbP_\alpha$-name by clause (b) and $\Vdash_{\bbP_\alpha}
``(I_\alpha[\bold x_\varepsilon],f_\alpha[\bold x_\varepsilon],
\name{\bbQ}_\alpha[\bold
x_\varepsilon]) \le^{\wk}_{\bold K} (I_\alpha,f_\alpha,\name{\bbQ}_\alpha)$
\sn
\item[$(h)$]   $(\alpha) \quad$ if $j < i$ then
$\Vdash_{\bbP_{\beta_i}} ``\name{\bar D}_{\beta_j} \le 
\name{\bar D}_{\beta_i}"$
\sn
\item[${{}}$]   $(\beta) \quad$ if $i$ is a limit ordinal and $t \in
\cT$ \then \, $\Vdash_{\bbP_{\beta_i}} ``\name D_{\beta_i,t} =$

\hskip25pt  $\cup\{\name D_{\beta_j,t}:j<i\}$
\sn
\item[$(i)$]   $(\alpha) \quad \langle \beta_j:j \le i\rangle$ is
increasing continuous
\sn
\item[${{}}$]  $(\beta) \quad$ if $i=0$ then $\beta_i=0$
\sn
\item[${{}}$]  $(\gamma) \quad$ if $i = j+1$ then
\sn
\item[${{}}$]  \hskip20pt $\bullet_1 \quad \beta_i \in S \cap E$
\sn
\item[${{}}$]  \hskip20pt $\bullet_2 \quad$ if $\gamma \in
[\beta_i,\lambda^+] \wedge \gamma \in (S \cap E) \cup \{\lambda^+\}$ 
and $t \in \cT$, \then \, 

\hskip40pt $\Vdash_{\bbR_{\beta_i,\gamma}[\bbP_{\beta_i},\bar{\bold x}]} 
``\emptyset \notin \fil(\cup\{\name D_{\gamma,t}
[\bold x_\varepsilon]:\varepsilon < \delta\} \cup \name D_{\beta_i,t})"$
\sn
\item[${{}}$]  \hskip20pt $\bullet_3 \quad$ if $\beta_j \in S \cap E$
then clause (g) of Definition \ref{c49} holds
\sn
\item[${{}}$]  \hskip20pt $\bullet_4 \quad$ if $\beta_j \in S_0$ and
$C^*_{\beta_j} \subseteq \{\beta_\iota:\iota < j\}$ \then \,
$\name{\bbQ}_{\beta_j} = \name{\bbQ}_{\fil(\emptyset)}$, and so

\hskip40pt  the relevant case of 
clause $(h)(\beta)$ of Definition \ref{c49} holds
\sn
\item[${{}}$]  $(\delta) \quad$ if $i$ is a limit ordinal, $\gamma \in
(\beta_i,\lambda^+] \wedge \gamma \in (S \cap E) \cup \{\lambda^+\}$ 

\hskip25pt and $t \in \cT$ \then \, $\Vdash_{\bbR_{\beta_i,\gamma}
[\bbP_{\beta_i},\bar{\bold x}]} ``\emptyset \notin
\fil(\cup\{\name D_{\gamma,t}[\bold x_\varepsilon]:\varepsilon <
\delta\}$

\hskip30pt $\cup \bigcup\{D_{\alpha,t}:\alpha < \beta_i\})"$.
\end{enumerate}
\end{enumerate}
\mn
Note that as $\name{\bar D}_\alpha$ (when $(\exists j \le
i)(\alpha = \beta_j \in S \vee j=0))$ is an 
ultra $\cT$-filter system, we do not have to
bother proving $A \in (\name D^+_{\alpha,s}[\name G_{\bbP_\alpha}])
\Rightarrow A \in (D^+_{\beta,s}[\name G_{\bbP_\beta}])$ (when $\alpha
< \beta$ are from $\{\beta_j:j \le i,\beta_j \in S\})$.

Also
\mn
\begin{enumerate}
\item[$(*)_1$]  if $t \in \cT,\varepsilon < \delta,\beta \le \beta_i$
and $\beta \in S \cap E_{\bold x_\varepsilon}$ then
$\Vdash_{\bbP_\beta} ``\name D^{\bold x_\varepsilon}_{\beta,t}
\subseteq \name D_{\beta_i,t}"$.
\end{enumerate}
\mn
[Why?  This follows from clause (i) of $\boxdot$.]

Let us carry the induction, this clearly suffices.
\bigskip

\noindent
\underline{Case 1}:  $i = 0$.

Trivial.
\bigskip

\noindent
\underline{Case 2}:  $i$ is a limit ordinal.

Let $\beta = \beta_i$ be $\cup\{\beta_j:j < i\}$, clearly $\langle
\beta_j:j \le i \rangle$ is increasing continuous and $\beta_i \in
E$.  Below $\varepsilon$ vary on $\delta$.

Let $\bbP_\beta = \cup\{\bbP_\alpha:\alpha < \delta\}$ and
$f_\beta = \cup\{f_\alpha:\alpha < \beta\}$ and from $\boxplus_0$
recall $I_{< \beta} =
\cup\{I_\alpha:\alpha < \beta\}$.  Clearly $\bbP_\beta
\in \bold K$ as witnessed by $(I_{< \beta},f_\beta)$ and $\alpha <
\beta \Rightarrow \bbP_\alpha \lessdot \bbP_\beta$.  Note that
$\bbP_\beta$ satisfies the c.c.c. as $\langle \bbP_\alpha:\alpha <
\beta\rangle$ is $\lessdot$-increasing continuous and the induction
hypothesis; alternatively using $f_\alpha$.

Now
\mn
\begin{enumerate}
\item[$(*)_2$]  $\bbP_\beta[\bold x_\varepsilon] \lessdot \bbP_\beta$
for $\varepsilon < \delta$; hence
$\bbR_{\beta,\gamma}[\bbP_\beta,\bar{\bold x}]$ is well defined for
$\gamma \in [\beta,\lambda^+]$.
\end{enumerate}
\mn
[Why?  Again we shall use \ref{z4}(3).

First, $\bbP_\beta[\bold x_\varepsilon] = \cup\{\bbP_{\beta_j}[\bold
x_\varepsilon]:j<i\}$ but $j<i \Rightarrow \bbP_{\beta_j}[\bold
x_\varepsilon] \subseteq \bbP_{\beta_j} \subseteq \bbP_\beta$
so clearly $\bbP_\beta[\bold x_\varepsilon] \subseteq \bbP_\beta$.

Second, $\bbP_\beta[\bold x_\varepsilon] \subseteq_{\ic} \bbP_\beta$,
because if $p,q \in \bbP_\beta[\bold x_\varepsilon]$ are incompatible
in $\bbP_\beta[\bold x_\varepsilon]$ then
for some $j<i$ we have $p,q \in \bbP_{\beta_j}[\bold x_\varepsilon]$
hence $p,q$ are incompatible in $\bbP_{\beta_j}[\bold x_\varepsilon]$,
so as $\bbP_{\beta_j}[\bold x_\varepsilon] \subseteq_{\ic}
\bbP_{\beta_j}$ they are incompatible in $\bbP_{\beta_j}$, but 
$\bbP_{\beta_j} \lessdot \bbP_\beta$ so they are incompatible in
$\bbP_\beta$ as required.

Third, if $q \in \bbP_\beta$ then for some $\alpha(0) < \beta$ we have $q
\in \bbP_{\alpha(0)}$ and so there is $p \in \bbP_{\alpha(0)}
[\bold x_\varepsilon]$ such that $p \le p' \in \bbP_{\alpha(0)}[\bold
x_\varepsilon] \Rightarrow p',q$ are compatible in $\bbP_{\alpha(0)}$.
 So it suffices to prove $p \le p' \in \bbP_\beta[\bold x_\varepsilon]
\Rightarrow p',q$ are compatible in $\bbP_\beta$, so fix such $p'$.
As $\beta$ is a limit ordinal, $\bbP_\beta = \cup\{\bbP_\alpha:\alpha <
\beta\}$ hence there is $\alpha(1)$ such that $\alpha(0) \le \alpha(1)
< \beta$ and $p' \in \bbP_{\alpha(1)}[\bold x_\varepsilon]$.  Now $p''
:= p' \rest \alpha(0)$ is well defined and belong to
$\bbP_{\alpha(0)}[\bold x_\varepsilon]$ and is above $p$, so by the
choice of $p$ there is a common upper bound $q^+ \in \bbP_{\alpha(0)}$
of $q$ and $p''$.  As $\langle \bbP_\alpha,\name{\bbQ}_\alpha:\alpha <
\beta\rangle$ is FS iteration, $q^+ \in \bbP_{\alpha(0)},p' \in
\bbP_{\alpha(1)}[\bold x_\varepsilon] \lessdot \bbP_{\alpha(1)}$ and
$p' \rest \alpha(0) \le q^+$, clearly there is a common upper bound
$r \in \bbP_{\alpha(1)} \lessdot \bbP_\beta$ of $p',q^+$ so $r$
exemplifies $p',q$ are compatible in $\bbP_\beta$.  So we have
finished proving $(*)_2$.]

Let $\name D'_{\beta,t} = \cup\{\name D_{\alpha,t}:\alpha = \beta_j$
 for some $j<i$ so $\alpha < \beta\}$.  Clearly $s
\le_{\cT} t \Rightarrow \name D'_{\beta,s} \subseteq \name
D'_{\beta,t}$ so the main point is to prove not just
$\Vdash_{\bbP_\beta} ``\emptyset
\notin \fil(\name D'_{\beta,t})"$, but that moreover $\gamma \in
[\beta,\lambda^+] \wedge \gamma \in (S \cap E) \cup \{\lambda^+\}
\Rightarrow \Vdash_{\bbR_{\beta,\gamma}[\bbP_\beta,\bar{\bold x}]} ``\emptyset
\notin \fil(D'_{\beta,\gamma,t})"$ where $D'_{\beta,\gamma,t} =
\cup\{\name D_{\gamma,t}[\bold x_\varepsilon]:\varepsilon < \delta\}
\cup D'_{\alpha,t} = \cup\{\name D_{\gamma,t}[\bold x_\varepsilon]:
\varepsilon < \delta\} \cup \bigcup\{\name D_{\alpha,t}:\alpha =
\beta_j$ for some $j<i\}$.
Fixing such $\gamma$, again as 
$\langle \name D^{\bold x_\varepsilon}_{\gamma,t}:
\varepsilon < \delta\rangle$ is increasing
and $\langle \name D_{\alpha,t}:\alpha = \beta_j$ for some
$j<i \rangle$ is increasing, it
suffice to prove $\Vdash_{\bbR_{\beta,\gamma}[\bbP_\beta,\bar{\bold x}]} 
``\emptyset \in \fil(\name D^{\bold x_\varepsilon}_{\gamma,t} 
\cup \name D_{\alpha,t})"$, for any $\varepsilon < \delta$ and $\alpha
 = \beta_j,j < i$.  For this it suffices to prove:
\mn
\begin{enumerate}
\item[$(*)_3$]  if (A) then (B) where
\begin{enumerate}
\item[$(A)$]  $(a) \quad p = (p_1,p_2) \in
\bbR_{\beta,\gamma}[\bbP_\beta,\bar{\bold x}]$
\sn
\item[${{}}$]  $(b) \quad t \in \cT$
\sn
\item[${{}}$]  $(c) \quad \alpha = \beta_j < \beta$ and $\name A \in \name
D_{\alpha,t}$ a $\bbP_\alpha$-name of a subset of $\bbN$
\sn
\item[${{}}$]  $(d) \quad \varepsilon < \delta$ and $\name B \in \name
D^{\bold x_\varepsilon}_{\gamma,t}$ a 
$\bbP^{\bold x_\varepsilon}_\gamma$-name of a subset of $\bbN$
\sn
\item[${{}}$]  $(e) \quad n_* \in \bbN$
\sn
\item[$(B)$]  $p \Vdash_{\bbR_{\beta,\gamma}[\bbP_\beta,\bar{\bold x}]} 
``\name A \cap \name B \nsubseteq [0,n_*)"$.
\end{enumerate}
\end{enumerate}
\bigskip

\noindent
\underline{Proof of $(*)_3$}:

Let $(\varepsilon_0,p_0)$ be a witness for $(p_1,p_2) \in
\bbR_{\beta,\gamma}[\bbP,\bar{\bold x}]$; as we can increase
$\varepsilon_0$, by $\boxplus_3$, and we can increase $\varepsilon$,
\wilog \, $\varepsilon_0 = \varepsilon$.

Without loss of generality $p_0,p_2 \in \bbP_\alpha$, as we can
increase $\alpha$, moreover as $\iota < i \Rightarrow 
\beta_{\iota +1} \in E \cap S$, similarly \wilog \,
$\alpha \in S \cap E$.  Let $p^*_2 \in \bbP_\alpha$ be a common upper
bound of $p_0,p_2$.  
We define a $\bbP^{\bold x_\varepsilon}_\alpha$-name $\name A'$ by:
\mn
\begin{enumerate}
\item[$(*)_{3.1}$]  if $\bold G \subseteq \bbP^{\bold
x_\varepsilon}_\alpha$ is generic over $\bold V$ then $\name A'[\bold
G] = \{n$: some $q \in \bbP_\alpha/\bold G$ forces 
$n \in \name A$ and if $p^*_2 \in \bbP_\alpha/\bold G$ then
$\bbP_\alpha \models ``p^*_2 \le q"\}$.
\end{enumerate}
\mn
Easily
\mn
\begin{enumerate}
\item[$(*)_{3.2}$]  $\name A'$ is a $\bbP^{\bold
x_\varepsilon}_\alpha$-name of a subset of $\bbN$
\sn
\item[$(*)_{3.3}$]  $\Vdash_{\bbP_\alpha} ``\name A \subseteq \name A'$".
\end{enumerate}
\mn
As $\bold x_\varepsilon \in \bold Q$ and $\alpha \in S \cap E
\subseteq S \cap E_{\bold x_\varepsilon}$ and $\bbP^{\bold
  x_\varepsilon}_\alpha \lessdot \bbP_\alpha$ and
$\Vdash_{\bbP_\alpha} ``\name D^{\bold x_\varepsilon}_{\alpha,t}
\subseteq \name D_{\alpha,t}"$, it follows that
\mn
\begin{enumerate}
\item[$(*)_{3.4}$]  $\Vdash_{\bbP_\alpha[\bold x_\varepsilon]} ``\name
A' \in \name D^{\bold x_\varepsilon}_{\alpha,t}$".
\end{enumerate}
\mn
But $\bbP_\alpha[\bold x_\varepsilon] \lessdot \bbP_\gamma[\bold
x_\varepsilon]$ hence, recalling $(A)(\alpha)$ of $(*)_3$:
\mn
\begin{enumerate}
\item[$(*)_{3.5}$]  $\Vdash_{\bbP_\gamma[\bold x_\varepsilon]} ``\name A' \in 
\name D_{\gamma,t}[\bold x_\varepsilon]$ 
and $\name B \in \name D_{\gamma,t}[\bold x_\varepsilon]$"
\end{enumerate}
\mn
hence
\mn
\begin{enumerate}
\item[$(*)_{3.6}$]  $\Vdash_{\bbP_\gamma[\bold x_\varepsilon]} 
``\name A' \cap \name B \in \name D_{\gamma,t}[\bold x_\varepsilon]$".
\end{enumerate}
\mn
Let $p'_0 \in \bbP_\beta[\bold x_\varepsilon]$ be such that $p'_0 \le
p'' \in \bbP_\beta[\bold x_\varepsilon] \Rightarrow p'',p^*_2$ are
compatible and \wilog \, $p_0 \le p'_0$.  Let $p_3 \in 
\bbP_\gamma[\bold x_\varepsilon]$ be above $p_1$ and
$p'_0$; by $(*)_{3.6}$ there are $q_1$ and $n$ such that: $p_3 \le q_1
\in \bbP_\gamma[\bold x_\varepsilon]$ and $n \ge n_*$ and $q_1
\Vdash_{\bbP_\gamma[\bold x_\varepsilon]} ``n \in \name A' \cap \name B"$.

Let $q_0 = q_1 \rest \beta$, it belongs to to $\bbP_\beta[\bold
x_\varepsilon]$; clearly $q_0 \le q \in \bbP_\beta[\bold
x_\varepsilon] \Rightarrow p^*_2,q$ are compatible in $\bbP_\gamma[\bold
x_\varepsilon]$.  Also clearly $p'_0 \le q_0 \in \bbP_\alpha[\bold
x_\varepsilon]$ so there is $r_1$ such that $q_0 \le r_1 \in
\bbP_\alpha[\bold x_\varepsilon]$ and $r_1$ forces a truth value to
``$n \in \name A'$" so as $r_1$ is compatible with $q_1$, necessarily
$r_1 \Vdash ``n \in \name A'"$.  So $p_0 \le q_0 \le r_1 \in
\bbP_\beta[\bold x_\varepsilon]$. 

By  the definition of $\name A'$ and the choice of $p_0$, there is $q_2 \in
\bbP_\gamma[\bold x_\varepsilon]$ such that:
\mn
\begin{enumerate}
\item[$(*)_{3.7}$]  $(a) \quad \bbP_\alpha \models
``p^*_2 \le q_2$ and $q_0 \le r_1 \le q_2"$
\sn
\item[${{}}$]  $(b) \quad q_2 \Vdash_{\bbP_\alpha} ``n \in \name A"$.
\end{enumerate}
\mn
Let $\alpha(1) < \beta$ be $\ge \alpha$ such that $q_1 \rest \beta \in
\bbP_{\alpha(1)}[\bold x_\varepsilon]$; as $(q_1 \rest \beta) \rest
\alpha = q_0 \le r_1 \le q_2$ and as $\langle
\bbP_\gamma,\name{\bbQ}_\gamma:\gamma < \beta \rangle$ is a FS
iteration, clearly $q_1 \rest \beta,q_2$ are compatible in
$\bbP_{\alpha(1)}$ and let $q_4 \in \bbP_{\alpha(1)}$ be a common
upper bound of $(q_1 \rest \beta),q_2$.  Let $q'_0 \in
\bbP_{\alpha(1)}[\bold x_\varepsilon]$ be such that $q'_0 \le q \in
\bbP_{\alpha(1)}[\bold x_\varepsilon] \Rightarrow q,q_4$ are
compatible in $\bbP_{\alpha(1)}$, 
so as $(q_1 \rest \beta) \le q_4$, \wilog \, $(q_1 \rest
\beta) \le q'_0$.
\mn
\begin{enumerate}
\item[$(*)_{3.8}$]  $q'_0 \in \bbP_\beta[\bold x_\varepsilon]$ and
$(\varepsilon,q'_0)$ witness $(q_1,q_4) \in
\bbR_{\beta,\gamma}[\bbP_\beta,\bar{\bold x}]$.
\end{enumerate}
\mn
[Why?  As $\bold x_\varepsilon \in \bold Q$ and $q_1 \rest \beta = q_0
\le q'_0$  clearly $q'_0
\Vdash_{\bbP_\beta[\bold x_\varepsilon]} ``q_1 \in \bbP_\gamma[\bold
x_\varepsilon]/\name G_{\bbP_\beta[\bold x_\varepsilon]}"$.

For proving $q'_0 \Vdash_{\bbP_\beta[\bold x_\varepsilon]} ``q_4 \in
\bbP_\beta/\name G_{\bbP_\beta[\bold x_\varepsilon]}"$ recall the
choice of $q'_0$.]
\mn
\begin{enumerate}
\item[$(*)_{3.9}$]  $(q_1,q_4) \Vdash_{\bbR_{\beta,\gamma}
[\bbP_\beta,\bar{\bold x}]} ``n \in \name A \cap \name B \backslash
[0,n_*)"$.
\end{enumerate}
\mn
[Why?  First, $q_1 \Vdash_{\bbP_\gamma[\bold x_\varepsilon]} ``n \in
\name B"$ by the choice of $q_1$ hence $(q_1,q_4)
\Vdash_{\bbR_{\beta,\gamma}[\bbP_\beta,\bar{\bold x}]} ``n \in \name B"$
recalling $\bbP_\gamma[\bold x_\varepsilon]
\lessdot \bbR_{\beta,\gamma}[\bbP_\beta,\bar{\bold x}]$ by
$\boxplus_7$.

Second, $q_4 \Vdash_{\bbP_\beta} ``n \in \name A"$ because $q_2
\Vdash_{\bbP_\alpha} ``n \in \name A"$ and $q_2 \le q_4,\bbP_\alpha \lessdot
\bbP_\beta$ and so $(q_1,q_4)
\Vdash_{\bbR_{\beta,\gamma}[\bbP_\beta,\bar{\bold x}]} ``n \in \name
A"$ because $\bbP_\beta \lessdot
\bbR_{\beta,\gamma}[\bbP_\beta,\bar{\bold x}]$ by $\boxplus_6$.

Third, $n \ge n_*$ recalling the choice of $n$.  So $(*)_{3.9}$ holds.]

Together we have proved $(*)_3$.

Lastly, clearly $\beta_i \in E$ and let $\name{\bar
D}_\beta = \name{\bar D}'_\beta$.  If $\beta = \beta_i \notin S$ we are
done.  So assume $\beta \in S$; by the induction hypothesis
$\alpha = \beta_j < \beta \Rightarrow
\Vdash_{\bbP_{\beta_{j+1}}} ``\name{\bar D}_{\beta_{j +1}}$ is ultra
$\cT$-filter system", and $\name{\bar D}_\alpha$ increases with
$\alpha$, also necessarily $\cf(\beta) = \lambda$ hence $\Vdash_{\bbP_\beta}
``\langle \cup\{\name D_{\alpha,t}:\alpha < \beta\}:t \in \cT\rangle$
is ultra hence $\name{\bar D}'$ is ultra so we are done.
\bigskip

\noindent
\underline{Case 3}:  $i = j +1,\beta_j \notin S \cup S_0$.

Let $\gamma \in (\beta_j,\lambda^+]$ and $\bbR = 
\bbR_{\beta_j,\gamma}[\bbP_{\beta_j},\bar{\bold x}]$, recalling
$\boxplus_5$ we know $\bbR$ satisfies the c.c.c., 
by $\boxplus_6$ we know $\bbP_{\beta_j} \lessdot \bbR$ and by
$\boxplus_7$ we know $\varepsilon <
\delta \Rightarrow \bbP_\gamma[\bold x_\varepsilon] \lessdot
\bbR$.  For $t \in \cT$, let 
$\name D'_{\beta_j,\gamma,t} = \cup\{\name D_{\gamma,t}[\bold
x_\varepsilon]:\varepsilon < \delta\} \cup \name D_{\beta_j,t}$, noting
$\name D_{\beta_j,t} = \bigcup\{\name D_{\beta_\iota,t}:
\iota \le j\}$, so by the
induction hypothesis, $\Vdash_{\bbR} ``\emptyset \notin \fil(\name
D'_{\beta_j,\gamma,t})"$ so $\name{\bar D}'_{\beta_j,\gamma,t} =
\langle \name D'_{\beta_j,\gamma,t}:t \in \cT\rangle$ is a
$\bbR_{\beta_j,\gamma}[\bbP_{\beta_j}]$-name of a $\cT$-filter system.
Hence there is $\name{\bar D}''_{\beta_i,\gamma} = \langle
D''_{\beta_i,\gamma,t}:t \in \cT\rangle$, a $\bbR$-name of an ultra
$\cT$-filter system above $\bar D'_{\beta_i,\gamma}$, \wilog \,
$D''_{\beta_i,\gamma,t} = \fil(D''_{\beta_i,\gamma,t})$ for $t \in
\cT$.  In particular
this holds for $\gamma = \lambda^+$ hence $E^*_i$ is a club of
$\lambda^+$ where
\mn
\begin{enumerate}
\item[$(*)_4$]  $E^*_i = \{\gamma < \lambda^+:\gamma$ is a limit
ordinal from $E$ and if $\xi < \gamma$ then $\langle D''_{\beta_i,\lambda^+,t}
\cap \cP(\bbN)^{\bold V[\bbP_\xi]}:t \in \cT\rangle$ is a
$\bbR_{\beta_j,\xi_1}[\bbP_{\beta_j},\bar{\bold x}]$-name for 
some $\xi_1 < \gamma\}$ .
\end{enumerate}
\mn
So we can choose $\beta_i = \beta(i) \in E^*_i \cap E \cap S
\backslash(\beta_j+1)$.

Let $\bbP_{\beta_i} =
\bbR_{\beta_j,\beta_i}[\bbP_{\beta_j},\bar{\bold x}]$ and similarly
$\bbP_\alpha = \bbR_{\beta_j,\alpha}[\bbP_{\beta_j},\bar{\bold x}]$
for $\alpha \in (\beta_j,\beta_i)$  and $\name{\bar D}_{\beta_i} =
\langle D''_{\beta_i,\lambda^+,t} \cap \cP(\bbN)^{\bold
V[\bbP_{\beta(i)}]}:t \in \cT\rangle$.

Also the choice of $\name{\bbQ}_\alpha,\name g_\alpha$ for $\alpha \in
[\beta_j,\beta_i)$ is dictated by clause (g) of $\boxdot$ hence also
of $f_\alpha$ and it is easy to check that all the clauses in the
induction hypothesis are satisfied.
\bigskip

\noindent
\underline{Case 4}:   $i = j+1,\beta_j \in S$.

So $\Vdash_{\bbP_{\beta_j}} ``\name{\bar D}_{\beta_j}$ is an ultra
$\bbP_{\beta_j}$-filter system".  Let $\beta = \beta_j$.

Let $\name{\bbQ}_\beta = \name{\bbQ}_{\name D_\beta},
\bbP_{\beta +1} = \bbP_{\beta} * \bbQ_{\name D_{\beta}}$.  
By Claim \ref{cn.28}, $\bbP_{\beta +1}[\bold x_\varepsilon] =
\bbP_\beta[\bold x_\varepsilon] * \name{\bbQ}_{\name D_\beta
[\bold x_\varepsilon]} \lessdot \bbP_\beta * \bbQ_{\name D_\beta} =
\bbP_{\beta_j+1}$ for $\varepsilon < \delta$.  
So $\bbR_{\beta +1,\gamma}[\bbP_{\beta +1},\bar{\bold x}]$ is well
defined for $\gamma \in [\beta +1,\lambda^+]$.

For $t \in \cT$ let $\name D'_{\beta +1,s}$ be the dual of 
$\name \id_{\name{\bold d}_{t(\beta),s}}[\bbP_\beta]$, a 
$\bbP_{\beta +1}$-name.
\mn
\begin{enumerate}
\item[$(*)_5$]
$\Vdash_{\bbR_{\beta +1,\gamma}[\bbP_{\beta +1},\bar{\bold x}]}
``\emptyset \notin \fil(\cup\{\name D_{\gamma,s}[\bold x_\varepsilon]:
\varepsilon < \delta\} \cup \name D'_{\beta,s}"$ for 
$\gamma \in (\beta,\lambda^+]$.
\end{enumerate}
\mn
Note that for $(\beta,\gamma)$ we know the parallel statements.  
\mn
\begin{enumerate}
\item[$(*)_6$]  convention: we write $(p_1,p_2,\name p_3) =
(p_1,(p_2,\name p_3))$ for members of $\bbR_{\beta
+1,\gamma}[\bbP_{\beta +1},\bar{\bold x}]$, where we treat
$\bbP_{\beta +1}$ as $\bbP_\beta * \name{\bbQ}_{\bar D_{\beta}}$,
so $p_2 \in \bbP_\beta$ and $\Vdash_{\bbP_\beta} ``\name p_3 \in
\name{\bbQ}_{\bar D_\beta}"$ and $\tr(\name p_3)$ is an object not just
a name.
\end{enumerate}
\mn
We need
\mn
\begin{enumerate}
\item[$(*)_7$]  if $(A)$ then $(B)$ where
\begin{enumerate}
\item[$(A)$]  $(a) \quad p = (p_1,p_2,\name p_3) \in \bbR_{\beta
+1,\gamma}[\bbP_{\beta +1},\bar{\bold x}]$
\sn
\item[${{}}$]  $(b) \quad t \in \cT$
\sn
\item[${{}}$]  $(c) \quad \name A$ is a $\bbP_{\beta +1}$-name of a member of
$\name D'_{\beta,t}$ that is, 

\hskip25pt $\Vdash_{\bbP_{\beta +1}} ``\name A$ is
$\name{\bold d}_{t(\beta),t}$-null"
\sn
\item[${{}}$]  $(d) \quad \varepsilon < \delta$ and $\name B$ is a
$\bbP_\gamma[\bold x_\varepsilon]$-name of a memer of $\name
D_{\gamma,t}[\bold x_\varepsilon]$
\sn
\item[${{}}$]  $(e) \quad n_* \in \bbN$
\sn
\item[$(B)$]   $p \Vdash_{\bbR_{\beta +1,\gamma}[\bbP_{\beta
+1},\bar{\bold x}]} ``\name A \cap \name B \nsubseteq [0,n_*)"$.
\end{enumerate}
\end{enumerate}
\mn
First note
\mn
\begin{enumerate}
\item[$(*)_{7.1}$]  $(a) \quad$ let $(\varepsilon,(p_0,p'_3))$ where
$(p_0,\name p'_3) \in \bbP_\beta[\bold
x_\varepsilon] * \bbQ_{\name{\bar D}_\beta[\bold x_\varepsilon]}$
witness 

\hskip25pt $p \in \bbR_{\beta +1,\gamma}[\bbP_{\beta +1},\bar{\bold x}]$
\sn
\item[${{}}$]  $(b) \quad$ let $q_2 \in \bbP_\beta$ be above $p_0,p_2$
\sn
\item[${{}}$]  $(c) \quad$ let $q_0 \in \bbP_\beta$ be such that $q_0
\le q' \in \bbP_\beta[\bold x_\varepsilon] \Rightarrow q_2,q'$ are
compatible
\sn
\item[${{}}$]  $(d) \quad$ let $\name B'$ be the following
$\bbP_{\beta +1}[\bold x_{\varepsilon +1}]$-name 

\hskip25pt $\{n$: there is $q \in \bbP_\gamma[\bold x_\varepsilon]/
\name{\bold G}_{\bbP_{\beta +1}[\bold x_\varepsilon]}$ forcing 
$n \in \name B$ above $p_1$ 

\hskip25pt when $p_1 \in \bbP_\gamma[\bold
x_\varepsilon]/\name{\bold G}_{\bbP_{\beta +1}[\bold x_\varepsilon]}\}$.
\end{enumerate}
\mn
Next consider
\mn
\begin{enumerate}
\item[$(*)_{7.2}$]  $\Vdash_{\bbP_{\beta +1}} ``\name A \cap \name B' \nsubseteq
[0,n_*)"$.
\end{enumerate}
\mn
Why is $(*)_{7.2}$ true?  Note that $\Vdash_{\bbP_{\beta +1}[\bold
x_\varepsilon]} ``\name B' \in (\name\id_{\bold d(t_\beta,s)})^+"$ by
clause (g) of Definition \ref{c49}, as $\Vdash_{\bbP_\gamma[\bold
x_\varepsilon]} ``\name B \in \name D_{\gamma,s}$ and $\name B'
\subseteq \name B"$.  Now apply Claim \ref{c37} for $\bbP_{\beta
+1}[\bold x_\varepsilon] = \bbP_\beta[\bold x_\varepsilon] *
\bbQ_{\name D_{\beta,t(\alpha)}[\bold x_\varepsilon]}$ and $\bbP_{\beta +1}
= \bbP_\beta * \bbQ_{\name D_{\beta,t(\alpha)}}$.  

Why is $(*)_{7.2}$ enough for proving $(*)_7$?  
As in the proof of Case 2, only much easier.
\bigskip

\noindent
\underline{Case 5}:  $i=j+1,\beta_j \in S_0$.

Let $\beta = \beta_j$; and let $\bbP_{\beta +1} = \bbP_\beta *
\name{\bbQ}_{\fil(\emptyset)}$ so $\name{\bbQ}_\beta =
\name{\bbQ}_{\fil(\emptyset)}$, and again $\bbP_\beta[\bold
x_\varepsilon] * \name{\bbQ}_{\fil(\emptyset)} \lessdot \bbP_\beta *
\name{\bbQ}_{\fil(\emptyset)}$ by \ref{cn.28}.  Clearly $\bbR_{\beta
+1,\gamma}[\bbP_{\beta +1},\bar{\bold x}]$ is well defined for $\gamma
\in [\beta +1,\lambda^+]$.  We let $D'_{\beta +1,t} = \cup\{\name
D_{\alpha,t}:\alpha \in S \cap E\} \cup\{\name u_{\beta,s,n}:
s \le_{\cT} t$ and $n \in \bbN\}$, a $\bbP_{\beta +1}$-name.

We have to prove the parallel of $(*)_5$, i.e.
\mn
\begin{enumerate}
\item[$(*)_8$]  $\Vdash_{\bbR_{\beta +1,\gamma}[\bbP_{\beta
+1},\bar{\bold x}]} ``\emptyset \notin \fil(\name D'_{\alpha,t})"$ for
$\gamma \in [\beta +1,\lambda^+]$ and $t \in \cT$.
\end{enumerate}
\mn
By \ref{c52} it suffices to prove
\mn
\begin{enumerate}
\item[$(*)_9$]  $\Vdash_{\bbR_{\beta +1,\gamma}[\bbP_{\beta +1},
\bar{\bold x}]} ``\emptyset \notin \fil(\{\name u_{\beta,s}:s
\le_{\cT} t\})$ for $t \in \cT$.
\end{enumerate}
\mn
Now it is like Case 4 only easier.
\end{PROOF}

\begin{claim}
\label{cn.77} 
If $\bold x \in \bold Q$ and 
$\theta \in \Theta_2$ \then \, we can find a pair $(\bold y,\bold j_*)$
such that
\mn
\begin{enumerate}
\item[$(a)$]   $\bold x \le_{\bold Q} \bold y$
\sn
\item[$(b)$]   $\bold j_*$ is an isomorphism from 
$(\bbP^{\bold x})^\theta/E_\theta$ onto $\bbP^{\bold y}$ 
extending $\bold j^{-1}_{**}$ where $\bold j_{**}$ is
the canonical embedding of $\bbP^{\bold x}_\kappa$ into 
$(\bbP^{\bold x}_\kappa)^\theta/E_\theta$
\sn
\item[$(c)$]   $\bold j_*$ maps 
$(\bbP^{\bold x}_\alpha)^\theta/E_\theta$ onto $\bbP^{\bold y}_\alpha$ 
for any $\alpha < \lambda^+$ satisfying $\cf(\alpha) \ne \theta$
\sn
\item[$(d)$]    note that $\bold j_*$ maps $\bold j_{**}
(\bbP^{\bold x}_\alpha)$ to a
$\lessdot$-subforcing of $\bbP^{\bold y}_\alpha$ for $\alpha \le
\lambda^+$ satisfying $\cf(\alpha) \ne \theta$.
\end{enumerate}
\end{claim}

\noindent
Before proving \ref{cn.77} recall
\begin{definition}
\label{p.1b}  
1) For a c.c.c. forcing notion $\bbP$ and $\bbP$-name $\name A$ of 
a subset of $\bbN$ we say that $\bold p = \langle (p_{n,m},
\bold t_{n,m}):m,n < \omega \rangle$ represents $\name A$ \when \,:
\mn
\begin{enumerate}
\item[$(a)$]   $p_{n,m} \in \bbP$ and $\bold t_{n,m}$ is a truth value
\sn
\item[$(b)$]   for each $n,\langle p_{n,m}:m < \omega \rangle$ is a
maximal antichain of $\bbP$
\sn
\item[$(c)$]   for $n,m < \omega$ we have $p_{n,m} \Vdash_{\bbP}
``n \in \name A$ \Iff \,  $\bold t_{n,m}"$.
\end{enumerate}
\mn
2) For $\bold p$ as in part (1) let $\name A_{\bold p}$ be 
the canonical $\bbP$-name represented by $\bold p$.
\end{definition}

\begin{fact}
\label{cn.81}  
1) If $\bbP$ is a c.c.c. forcing notion and
$\name A$ is a $\bbP$-name of a subset of $\bbN$ \then \, 
some $\langle (p_{n,m},\bold t_{n,m}):n,m < \omega \rangle$
represents $\name A$. 

\noindent
2) If $\bbP$ is a c.c.c. forcing notion and $\name A',\name A''$ 
are $\bbP$-names of subsets of $\omega$, both represented by 
$\langle (p_{n,m},\bold t_{n,m}):n,m < \omega
\rangle$ \then \, $\Vdash_{\bbP} ``\name A' = \name A''"$.

\noindent
3) For a sequence $\bar{\bold t} = \langle \bold t_{n,m}:n,m < \omega
\rangle$ of truth values, for some formula $\varphi = 
\varphi^0_{\bar{\bold t}}(\bar x) \in \bbL_{\aleph_1,\aleph_1}
(\tau),\tau = \{\le\}$ where $\bar x = \langle
x_{n,m}:n < \omega \rangle$ we have: for every c.c.c. forcing notion
$\bbP$ and $p_{n,m} \in \bbP\, (n,m < \omega)$ we have:
\mn
\begin{enumerate}
\item[$\circledast$]   $\bbP \models 
``\varphi(\langle p_{n,m}:n,m < \omega \rangle)$ iff
$\langle (p_{_n,m},\bold t_{n,m}):n,m < \omega \rangle$ represents a
$\bbP$-name of a  non-empty subset of $\omega"$.
\end{enumerate}
\mn
4) For $k < \omega$, sequences $\bar{\bold t}^\ell = \langle
\bold t^\ell_{n,m}:n,m < \omega \rangle$ of truth values for $\ell \le k$
for some $\bbL_{\aleph_1,\aleph_1}(\tau)$-formula 
$\varphi = \varphi^k_{\bold t^0,\dotsc,\bold t^k}(y,\bar
x^0,\dotsc,\bar x^k)$ where $\bar x^\ell = \langle x^\ell_{n,m}:n,m <
\omega \rangle$ we have:
\mn
\begin{enumerate}
\item[$\circledast$]   for every 
$q,p^\ell_{n,m} \in \bbP \, (n,m < \omega,\ell \le k),\bbP$ 
a c.c.c. forcing notion we have: 

$\bbP \models \varphi[q,\langle p^0_{n,m}:n,m < \omega \rangle,
\langle p^1_{n,m}:n,m < \omega
\rangle,\ldots,\langle p^k_{n,m}:n,m < \omega \rangle]$ \Iff \, $\langle
(p^\ell_{n,m},\bold t^\ell_{n,m}):n,m < \omega \rangle$ represents a
$\Bbb P$-name of a subset of $\omega$ which we call $\name A_\ell$, 
for $\ell \le k$ and $q \Vdash_{\bbP}
``\name A_k$ and $\bbN \backslash \name A_k$ do not almost
include $\name A_0 \cap \name A_1 \cap \ldots \cap \name A_{k-1}"$.
\end{enumerate}
\end{fact}

\begin{PROOF}{\ref{cn.81}}
Easy.
\end{PROOF}

\begin{remark}
\label{cn.83}
In \ref{cn.81} we can treat any other relevant properties of such
$\bbP$-names.
\end{remark}

\begin{PROOF}{\ref{cn.77}}
\underline{Proof of \ref{cn.77}}  

Let $\chi$ be large enough, $\bold x \in {\cH}(\chi)$ and ${\gB} =
({\cH}(\chi),\in)^\theta/E_\theta$ and let 
$\bold j$ the canonical embedding of $({\cH}(\chi),\in)$ into ${\gB}$.

We now define
\mn
\begin{enumerate}
\item[$\boxplus$]  $(a) \quad \bbP_\alpha$ is 
$(\bbP^{\bold j(\bold x)}_{\bold j(\alpha)})^{\gB}$ if $\alpha \le
\lambda^+,\cf(\alpha) \ne \theta$ and $\bbP_\alpha =
\cup\{\bbP_\beta:\beta < \alpha\}$ 

\hskip25pt if $\alpha < \lambda^+ \wedge \cf(\alpha) = \theta$
\sn
\item[${{}}$]  $(b) \quad I_{< \alpha} = \cup\{I^{\bold j(\bold x)}_{<
\bold j(\beta +1)})^{\gB}:\beta < \alpha\}$ and $E = E^{\bold y} =
E^{\bold x}$
\sn
\item[${{}}$]  $(c) \quad f_\alpha$, a function with domain
$\bbP_\alpha$ is defined by:
\sn
\begin{enumerate}
\item[${{}}$]  $(\alpha) \quad f_\alpha(p)$ is a function with domain $\{a:\gB
\models a \in \Dom(\bold j(f_\alpha(p))\}$
\sn
\item[${{}}$]  $(\beta) \quad f_\alpha(p)(a) = \bold
j^{-1}((f_\alpha(p)(a))^\gB)$
\end{enumerate}
\item[${{}}$]  $(d) \quad (I_\alpha,\name g_\alpha,\bbQ_\alpha)$ is
defined naturallly for $\alpha < \lambda^+$:
\sn
\begin{enumerate}
\item[${{}}$]  $(\alpha) \quad$ if $\cf(\alpha) \ne \theta$ as 
$(\bold j(\bbQ^{\bold x}_\alpha(p)))^{\gB}$
\sn
\item[${{}}$]  $(\beta) \quad$ if $\cf(\alpha) = \theta$, it is $\{p
\in \bbP_{\alpha +1}:\dom(p) \subseteq \bigcap\limits_{\beta < \alpha}
  [\bold j(\beta),\bold j(\alpha +1))^\gB\}$, 

\hskip25pt  etc.
\end{enumerate}
\item[${{}}$]  $(e) \quad E = E_{\bold x}$.
\end{enumerate}
\mn
We like to choose $\bbP^{\bold y}_\alpha = \bbP_\alpha$, a pedantic
objection is that $\bold j$ is not the identity, moreover $\bbP_\alpha
\subseteq \cH(\lambda^{++})$; so $\bbP^{\bold x}_\alpha 
\nsubseteq \bbP_\alpha$, by renaming we can overcome this.

Also for $\alpha \in E \cup \{\lambda^+\}$ and $t \in \cT$ the 
$\bbP^\bold y_\alpha$-name $\name D^{\bold y}_{\alpha,t}$ are naturally
defined such that
\mn
\begin{enumerate}
\item[$(*)$]  $\Vdash_{\bbP^{\bold y}_\alpha} 
``\name D^{\bold y}_{t,\alpha} = \{\name A_{\bold p}:\bold p$
represents some $\bbP^{\bold y}_{t,\alpha}$-name of subset of $\bbN$ and
$p \Vdash ``\bold p \in \bold j(\name D^{\bold
x}_{t,\alpha})"$ for some $p \in \name G_{\bbP^{\bold y}_\alpha}\}$.
\end{enumerate}
\mn
Almost all the desired 
properties hold by \L os theorem for $\bbL_{\aleph_1,\aleph_1}$ as in
\ref{cn.81}.   
A problem is to show clause $(d)(\alpha)$ of \ref{c49}, 
being ``ultra" which means
\mn
\begin{enumerate}
\item[$\odot$]   if $\partial \in \Theta_1,
s \in {\cT}_\partial,\alpha \in E \cap S$ then
$\Vdash_{\bbP^{\bold y}_\alpha}$ ``if $A \subseteq \bbN$ and $A
\ne \emptyset$ mod $\name D^{\bold y}_{\alpha,s}$ then for some $t$ we have $s
\le_{{\cT}_\partial} t$ and $A \in \name D^{\bold y}_{\alpha,t}"$.
\end{enumerate}
\mn
Toward this, as $\theta \in \Theta_2,\partial \in \Theta_1$ we have
$\theta \ne \partial$ hence
\mn
\begin{enumerate}
\item[$\boxdot$]   if $\eta$ be a generic branch of ${\cT}_\partial$
over $\bold V$ so $\eta$ is a subset of $\cT_\partial$ of order type
$\theta$ by $<_{\cT}$ \then \,
\begin{enumerate}
\item[$(\alpha)$]   $E_\theta$ is a $\theta$-complete ultrafilter
on $\theta$ even in $\bold V[\eta]$
\sn
\item[$(\beta)$]    $(\bbP_{\bold x})^\theta/E_\theta$ is the
same in $\bold V$ and $\bold V[\eta]$.
\end{enumerate}
\end{enumerate}
\mn
[Why?  The proof by the division to two cases:
\medskip

\noindent
\underline{First Case}:  $\theta < \partial$.  

The forcing $\cT_\partial$ adds to $\bold V$ no 
sequence of length $< \partial$ so obvious.
\medskip

\noindent
\underline{Second Case}:  $\theta > \partial$.

Note that $\bold j \restriction {\cT}_\partial$ is an isormorphism
from ${\cT}_\partial$ onto $(\bold j \restriction 
{\cT}_\partial)^{\gB}$ as $|{\cT}_\partial| < \theta$.]

So by $\boxdot$
\mn
\begin{enumerate}
\item[$\boxtimes$]   $\Vdash_{\cT_\partial}$ 
``$\{\name D^{\bold y}_{\alpha,\eta_t}:t \in \name \eta\}$ is an
ultrafilter on $\bbN$".
\end{enumerate}
\mn
This suffices for $\odot$ by \ref{cn.17} so we are done.
\end{PROOF}

\noindent
We lastly arrive to the desired conclusion.
\begin{conclusion}
\label{cn.91}  There is $\bbP$ such that (for our $\cT_*$ see
\ref{cn.42}(g), \ref{cn.47}):
\mn
\begin{enumerate}
\item[$(a)$]   $\bbP$ is a c.c.c. forcing notion 
of cardinality $\lambda$ and $\Vdash_{\bbP} ``2^{\aleph_0} = \lambda"$
\sn
\item[$(b)$]   ${\cT}_*$ has cardinality $\Pi\Theta_1 \le \lambda^+$, add no
new sequence of length $< \min(\Theta_1)$ of ordinals, 
collapse no cardinal, change no cofinality
\sn
\item[$(c)$]  $\bbP \times {\cT}_*$ has cardinality $\le \lambda
  + \Pi \Theta_1$, collapse no cardinality, change 
no confinality and forces $2^{\aleph_0} = \lambda$
\sn
\item[$(d)$]   in $\bold V^{\bbP \times {\cT}_*}$ we have $\Theta_1
\subseteq \Sp_\chi$, i.e. for every $\theta \in \Theta_1$ there
is a non-principal ultrafilter $D$ of character $\theta$
\sn
\item[$(e)$]   in $\bold V^{\bbP \times {\cT}_*}$ we have $\Theta_2
\cap \Sp_\chi = \emptyset$
\sn
\item[$(f)$]   $\bbP = \bbP^{\bold x}_{\delta(*)}$ for some 
$\bold x \in \bold Q$ and $\delta(*) \in E_{\bold x} \cap S$.
\end{enumerate}
\end{conclusion}

\begin{remark}
\label{c93}
1) So if $\sup(\Theta_1)$ is strongly Mahlo then $|\cT_*| =
   \sup(\Theta_1)$.

\noindent
2) Similarly in \ref{e.4} for $\Theta$.
\end{remark}

\begin{PROOF}{\ref{cn.91}}
We choose $\bold x_\varepsilon \in \bold Q$ by induction
on $\varepsilon \le \lambda$ such that
\mn
\begin{enumerate}
\item[$(*)$]   $(a) \quad \bold x_\varepsilon \in \bold Q$
\sn
\item[${{}}$]    $(b) \quad \zeta < \varepsilon \Rightarrow \bold x_\zeta
\le \bold x_\varepsilon$
\sn
\item[${{}}$]   $(c) \quad$ if $\varepsilon = \zeta +1$ and
cf$(\zeta) = \theta \in \Theta_2$ or cf$(\zeta) \notin \Theta_2 \wedge
\theta = \text{ min}(\Theta_2)$ then $\bold x_\varepsilon$

\hskip25pt  is gotten
from $\bold x_\zeta$ as $\bold y$ was gotten from $\bold x$ in
\ref{cn.77} using $E_\theta$
\sn
\item[${{}}$]  $(d) \quad \zeta < \varepsilon \Rightarrow
E_{\bold x_\varepsilon} \subseteq E_{\bold x_\zeta}$.
\end{enumerate}
\mn
For $\varepsilon=0$ use \ref{cn.51}, for $\varepsilon$ successor use
\ref{cn.77} and for $\varepsilon$ limit use \ref{cn.56}. 

Having carried the induction, let $\bold x = \bold x_\lambda$.  Let
$S'_0 = \{\delta \in S_0:C^*_\delta \subseteq E_{\bold x}\}$ so a
stationary subset of $\lambda^+$.
Let $E = \{\delta \in E_{\bold x}:\delta = \sup(\delta \cap S'_0)\}$.
Let $\delta(*) \in E$ be such that $\delta(*)$ has cofinality $\kappa$.
Let $\langle \alpha(\varepsilon):\varepsilon < \kappa\rangle$ be an
increasing sequence of members of $E_{\bold x}$ 
with limit $\delta(*)$ such that
$\varepsilon < \kappa \Rightarrow \alpha(\varepsilon +1) \in S'_0$.

Now letting $\bbP = \bbP^{\bold x}_{\delta(*)}$ recalling 
$\bbP^{\bold x}_{\delta(*)} = \cup\{\bbP^{\bold x_\varepsilon}_{\delta(*)}:
\varepsilon < \lambda\}$, it easily satisfies all the requirements but
we give some details.  We have $\Vdash_{\bbP} ``2^{\aleph_0} \ge
\lambda"$ and $|\bbP| \ge \lambda$
by the choice of $\bold x_0$ as $\bbP_1[\bold x_0] \lessdot \bbP$, 
see \ref{cn.51}; also $\bbP$ satisfies the
c.c.c. (see \ref{cn.53}(2)) and $\bbP$ has cardinality $\le \lambda$,
(see Definition \ref{c49}, clause (a)) hence $\Vdash_{\bbP}
``2^{\aleph_0} \le \lambda"$ recalling $\lambda =
\lambda^{\aleph_0}$.  So we have shown clause (a) of the
conclusion.  Clause (b) holds by the choice of $\cT_*$ (see end of
clause (g) of the
hypothesis \ref{cn.42}).  Now $|\bbP| = \lambda,|\cT| \le \Pi \Theta_1$
hence $|\bbP \times \cT_*| \le \lambda + \Pi \Theta_1$ 
and $\Vdash_{\cT_*} ``\bbP$ satisfies the c.c.c." by 
Hypothesis \ref{cn.42}(g); hence forcing with $\bbP
\times \cT_*$ collapse no cardinal which forcing with $\cT_*$ does not
collapse; but as $\theta \in \Theta_1 \Rightarrow \theta = \theta^{<
\theta}$ and the use of Easton support in the product
$\cT_*$, forcing with $\cT_*$ collapse no cardinal.  Similarly forcing
with $\bbP \times \cT_*$ changes no cofinality; together clause (c) of
\ref{cn.91} holds.

As for clause (d), as $\cT_*$ is a product, forcing with $\cT_*$ adds
$\name{\bar\eta} = \langle \name\eta_\theta:\theta \in
\Theta_1\rangle,\name\eta_\theta$ a $\theta$-branch of $\cT_\theta$ so in
$\bold V[\bar\eta]$ we have $\cup\{\name D^{\bold x_\lambda}_{\delta(*),t}:t
\in \name\eta_0\}$, which is a $\bbP$-name $\name D_\theta$ of an
ultrafilter on $\bbN$ by \ref{cn.17}(2), non-principal by \ref{cn.1}(2).  Now for each
$t \in \cT_\theta$, the filter $\name D^{\bold x_\lambda}_{\delta(*),t}$
is (forced to be) generated by the $\subseteq^*$-decreasing $\langle \name
u_{\alpha(\varepsilon +1),t,n}:\varepsilon < \kappa$ and $n \in \bbN\rangle$, 
in the sense that $\name u_{\alpha(\varepsilon +1),t,n+1} \subseteq
\name u_{\alpha(\varepsilon +1),t,n}$ and for $\zeta < \varepsilon$
for some $n_*$ we have $n_1 \in \bbN \wedge n_2 \in \bbN \backslash
n_* \Rightarrow \name u_{\alpha(\varepsilon +1),t,n_2} \subseteq^*
u_{\alpha(\zeta +1),t,n_1}$.  So $\name D_\theta$ is 
generated by $|\theta| + \kappa = \theta$
sets.  Now $\name\eta_\theta$ under $<_{\cT_\theta}$ has order type
$\theta$ and no $\name D^{\bold x_\lambda}_{\delta(*),t}$ is an
ultrafilter and it increases with $t$, so 
clearly $< \theta$ sets do not suffice.  Hence
$\Vdash_{\bbP \times \cT} ``\theta \in \Sp_\chi$ for every $\theta \in
\Theta_1"$, so clause (d) of \ref{cn.91} holds.

Lastly, concerning clause (e), assume that $(p,t) \in (\bbP
\times \cT_*)$ forces that ``$\name\cA \subseteq \cP(\bbN)$ generates a
non-principal ultrafilter $\name D$, of character $\theta,\theta =
|\name \cA|$ and $\theta \in \Theta_2"$.  As cf$(\lambda) > \theta$
and $\cT_* \equiv \cT_{\ge \theta} \times \cT_{< \theta}$ and $\cT_{\ge
\theta}$ is $\theta^+$-complete, min$(\Theta_1 \backslash \theta) >
\Pi(\Theta_1 \cap \theta) + \aleph_1$, \wilog \, $\name\cA$ is a
$(\bbP \times \cT_{< \theta})$-name.  As $\lambda \ge \cf(\lambda) > \theta \ge
\Pi(\Theta_1 \cap \theta)$ by \ref{cn.42}(1)(e) 
for some $\varepsilon < \lambda,\name\cA$ is a
$(\bbP^{\bold x_\varepsilon}_{\delta(*)} \times \cT_{< \theta})$-name.
As we can increase
$\alpha$ \wilog \, $\cf(\alpha) = \theta$.  Now apply \ref{cn.77}
recalling clause (c) of $(*)$.
\end{PROOF}
\newpage

\section {The $\aleph_n$'s and collapsing}

A drawback of \ref{cn.91} is that $\bold V$ and $\bold V^{\bbP}$
have the same cardinals while the cardinals missing from
$\Sp_\chi$ are ex-large cardinals so weakly inaccessible.  In
particular it gives no information on chaotic behaviour of $\Sp_\chi$
among the $\aleph_n$'s.  This is resolved to a large extent below.
However, here we do not improve the consistency strength, also we do
not deal here with successor of singulars but deal little with singulars.

So fulfilling the second promise from \S0 (the first was dealt with in
\S2, i.e. \ref{cn.91}) the main result of this section is
\begin{conclusion}
\label{e.16}  
1) If $u \subseteq \{1,2,\dotsc,n,\ldots\}$ and 
$n \ge 1 \Rightarrow n \in u \vee n +1
\in u$ and in $\bold V$ there are infinitely many measurable
cardinals, \then \, for some forcing notion $\bbP$ in $\bold
V^{\bbP}$ we have $\aleph_\omega \cap \Sp_\chi = \{\aleph_n:n
\in u\}$.

\noindent
2) Assume in $\bold V$ there are infinitely many compact cardinals.
Then in part (1) we can use any $u \subseteq [1,\omega)$.
\end{conclusion}

\begin{PROOF}{\ref{e.16}}
Straightforward from \ref{e.4}, \ref{e.8} below. 
\end{PROOF}

\begin{claim}
\label{e.4}  Assume G.C.H. for simplicity, Hypothesis
\ref{cn.42} and $\theta \in \Theta_2 \Rightarrow \theta > \sup(\theta
\cap \Theta)$ and $\cT_\theta$ is $\theta$-complete for $\theta \in
\Theta_1,\lambda = \cf(\lambda)$ for simplicity; 
let $\bold f$ be a function with
domain $\Theta_2$ such that $\theta > \bold f(\theta) > 
\sup(\Theta \cap \theta),\bold f(\theta) 
> \aleph_1$ is regular (so $\bold f(\theta)^{< \bold f(\theta)} 
= \bold f(\theta)$) and $\bold f(\theta) \notin
\Theta_2$ and let $\bbQ$ be the
product $\Pi\{\Levy(\bold f(\theta),< \theta):\theta \in \Theta_2\}$
with Easton support (recall {\rm Levy}$(\bold f(\theta),<\theta)$ is
collapsing each $\alpha \in [\bold f(\theta),\theta)$ to $\bold f(\theta)$ by
  approximation of cardinality $< \bold f(\theta))$.

Lastly, let $\bold x = \bold x_\lambda,\delta(*)$ be as in the proof of
\ref{cn.91}.
\Then \, $\bbP = \bbP^{\bold x}_{\delta(*)} \times {\cT}
\times \bbQ$ satisfies:
\mn
\begin{enumerate}
\item[$(a)$]   $\bbP$ is a forcing notion 
of cardinality $\lambda$ and $\Vdash_{\bbP} ``2^{\aleph_0} = \lambda"$
\sn
\item[$(b)$]   ${\cT}$ has cardinality $\le \Pi\Theta_1$, and as a
forcing notion adds no
new sequence of length $<$ {\rm min}$(\Theta_1)$ of ordinals, 
collapses no cardinal, changes no cofinality
\sn
\item[$(c)$]   $\bbP$ has cardinality $\le \lambda
  + \Pi \Theta_1$, really $\lambda + |\Pi \cT_*| + |\bbQ|$, collapses 
no cardinal except those in
$\cup\{(\bold f(\theta),\theta):\theta \in \Theta_2\}$, changes 
no confinality except that {\rm cf}$^{\bold V}(\delta) = (\bold
f(\theta),\theta) \Rightarrow \cf^{\bold V[\bbP]}(\delta) = \bold f(\theta)$.
\sn
\item[$(d)$]   In $\bold V^{\bbP}$ we have $\Theta_1
\subseteq \Sp_\chi$, i.e. for every $\theta \in \Theta_1$ there
is a non-principal ultrafilter $D$ of character $\theta$
\sn
\item[$(e)$]   in $\bold V^{\bbP}$ we have $\Theta_2
\cap \Sp_\chi = \emptyset$.
\end{enumerate}
\end{claim}

\begin{discussion}
\label{e5}
1) We may allow $\bold f(\theta) = \sup(\Theta \cap \theta)$ when
   $\sup(\Theta  \cap \theta) \notin \Theta_2$.

\noindent
2) We may like to have successive members of $\Theta_2$, see
\ref{e.8}; together with \ref{e5}(1) we get full answer for the $\aleph_n$'s.

\noindent
3) We may in \ref{e.4}, if $\lambda = \lambda^{< \kappa}$ demand
   $\Vdash_{\bbP} ``\MA_{< \kappa}"$, for this we need in the
   inductive choice of the $\bold x_\varepsilon$'s for $\varepsilon <
   \lambda$ another case; we do not get MA$_{\le \kappa}$ as
   $\cf(\delta(*)) = \kappa$.

\noindent
4) Similarly to part (3) in \ref{cn.19}, \ref{cn.91}, \ref{e.30}, \ref{e.16}.
\end{discussion}

\begin{PROOF}{\ref{e.4}}
First, clause (c), on when cardinals and cofinalities are preserved
should be clear.  Second, note that forcing by ${\cT}_* \times \bbQ$ adds 
no new $\omega$-sequence of members of $\bold V$ and even preserve
``$\bbP_{\bold x_\lambda}$ satisfies c.c.c." (and even ``satisfies the Knaster
condition" and even ``being locally $\aleph_1$-centered") all because
$\cT_*  \times \bbQ$ is $\aleph_1$-complete.  
So ${\cP}(\bbN)^{\bold V[\bbP]}$ and even 
$({}^\omega\Ord)^{\bold V[\bbP]}$ is the same 
as the one in $\bold V[\bbP^{\bold x}_{\delta(*)}]$.

Third, note that for every $\theta \in \Theta_1$, in 
$\bold V^{\cT_*}$ we have a 
$\bbP^{\bold x}_{\delta(*)}$-name $\name D_\theta$ of
an ultrafilter on $\bbN$ with $\chi(\name D_\theta) = \theta$, so
there is a set $\bbD_\theta$ of $\bbP^{\bold x}_{\delta(*)}$-names of reals of
cardinality $\theta$, or better a set of representations of such
names, (see Definition \ref{p.1b}), which generates $\name D_\theta$.

Now $\name D_\theta$ has the same properties 
in $\bold V^{\cT_* \times \bbQ}$ (see
``first" and ``second" above) so we have $\theta \in 
\Sp_\chi^{\bold V[\bbP]}$ so $\bold V^{\bbP} \models 
``\Theta_1 \subseteq \Sp_\chi"$.

Fourth, the main point, we would like to prove that 
$\Theta_2 \cap \Sp_\chi = \emptyset$ in $\bold V^{\bbP}$. 

So toward contradiction assume 
\mn
\begin{enumerate}
\item[$\odot_1$]   $\theta \in \Theta_2$ and $(p^*,r^*,q^*) \in 
\bbP$ forces ``$\name D$ is an ultrafilter on $\bbN$ with 
$\chi(\name D) = \theta"$.
\end{enumerate}
\mn
Let ${\bbQ}_{<\theta}$ be $\{p \in \bbQ:\dom(q) \subseteq
  \theta\}$ and similarly $\bbQ_{\le \theta},\bbQ_{> \theta}$ so essentially 
$\bbQ = \bbQ_{\le \theta} \times \bbQ_{> \theta}$ and
$\bbQ_{\le\theta} = \bbQ_{<\theta} \times \bbQ_\theta$ where
$\bbQ_\theta = \text{ Levy}(\bold f(\theta),<\theta)$.  
Similarly ${\cT}_{<\theta} = \{r \in {\cT}:\dom(r) \subseteq
\theta\}$, etc.

Now
\mn
\begin{enumerate}
\item[$(*)_1$]   $|{\cT}_{< \theta} \times \bbQ_{< \theta}| < \theta$.
\end{enumerate}
\mn
[Why?  Recalling $|\cT_{< \theta}| \le (\sup(\Theta_1 \cap \theta))^+
\le (\sup(\Theta \cap \theta))^+ \le \bold f(\theta)^+ < \theta$ 
by an assumption on $\bold f$ and $\bbQ_{< \theta} \le
\Pi\{\bbQ_\partial:\partial \in \Theta_2 \cap \theta\}$ has
cardinality $\le \sup(\Theta_2 \cap \theta)^+ \le \bold f(\theta)^+ < \theta$.]
\mn
\begin{enumerate}
\item[$(*)_2$]   there is a sequence
$\langle \name{\bold p}_\varepsilon:\varepsilon <
\theta\rangle,\name{\bold p}_\varepsilon$ a $(\cT_* \times \bbQ)$-name
of a $\bbP^{\bold x}_{\delta(*)}$-representation of a 
subset $\name A_\varepsilon$ of $\bbN$ such that $(p^*,r^*,q^*) \Vdash_{\bbP} 
``\{\name A_\varepsilon:\varepsilon < \theta\}$ generates $\name D$
and $\name A_n \cap [0,n) = \emptyset$ and $\chi(\name D) = \theta"$
\sn
\item[$(*)_3$]   \wilog \, $(p^*,r^*,q^*) \in \bbP' := 
\bbP^{\bold x}_{\delta(*)} \times {\cT}_{< \theta} \times \bbQ_{\le \theta}$
and $\name D$, moreover the sequence $\langle \name{\bold p}_\varepsilon:
\varepsilon < \theta\rangle$ are $\bbP'$-names.
\end{enumerate}
\mn
[Why?  Because, first, $\bbQ/\bbQ_{\le \theta}$ is 
$\theta^+$-complete as we are assuming $\sigma \in \Theta_2 \backslash
\theta^+ \Rightarrow \bold f(\sigma) > \theta$.  Second, recalling
$\theta \notin \Theta_1$ as $\Theta_1,\Theta_2$ are disjoint, forcing by
${\cT}_{\ge \theta} = \cT_{> \theta}$ adds no new sequence of length $\le
\theta$ of ordinals (by \ref{cn.42}) and even is $\theta^+$-complete
(by the claim assumptions).  Third, $\cT_{< \theta} \times
\bbQ_{\le \theta}$ has cardinality $\le \theta$.]
\mn
\begin{enumerate}
\item[$(*)_4$]   there are $\langle (r_\varepsilon,q_\varepsilon,
\bold q_\varepsilon,\name A'_\varepsilon):\varepsilon < 
\theta\rangle$ such that:
\begin{enumerate}
\item[$(a)$]   $r_\varepsilon \in {\cT}_{< \theta}$ and $q_\varepsilon 
\in \bbQ_{\le \theta}$
\sn
\item[$(b)$]   $\bold q_\varepsilon$ is a canonical representation of a
$\bbP^{\bold x}_{\delta(*)}$-name of a subset of $\bbN$
\sn
\item[$(c)$]   $(p^*,r_\varepsilon,q_\varepsilon)$ belongs to 
$\bbP^{\bold x}_{\delta(*)} \times {\cT}_{< \theta} \times \bbQ_{\le \theta}$, 
is above $(p^*,r^*,q^*)$ and forces that $\name A_{\bold q_\varepsilon},
\bbN \backslash \name A_{\bold q_\varepsilon}$ are $\ne \emptyset$ mod
fil$(\{\name A_\iota:\iota < \varepsilon\})$ and $\{\name
A_\iota:\iota < \varepsilon\}$ is included in this filter and the
condition also forces $\name{\bold p}_\varepsilon$ is $\bold q_\varepsilon$
\sn
\item[$(d)$]  $\name A'_\varepsilon$ is the $\bbP_{\bold x}$-name of a
subset of $\bbN$ represented by $\bold q_\varepsilon$
\sn
\item[$(e)$]   for technical reasons $\theta \in \dom(q^*_\varepsilon)$.
\end{enumerate}
\end{enumerate}
\mn
[Why?  As $(p^*,r^*,q^*)$ forces that $\{\name A_\iota:\iota <
  \theta\}$ generates $\name D$ but $\name A_\varepsilon \notin
\fil(\{\name A_\zeta:\zeta < \varepsilon\})$.]

Easily
\mn
\begin{enumerate}
\item[$(*)_5$]  there are representations $\bold q'_i(i < \theta)$ 
of $\bbP^{\bold x}_{\delta(*)}$-names $\name C_i$ such that
\begin{enumerate}
\item[$(a)$]  $(p^*,r^*,q^*) \Vdash_{\bbP'} ``\name{\bold p}_\varepsilon  
\in \{\bold q'_i:i <\theta\}"$ for every $\varepsilon < \theta$
\sn
\item[$(b)$]   $(p^*,r^*,q^*) \Vdash ``\{\name C_i:i < \theta\}$ 
includes $\{\name A_i:i < \theta\}$ and is closed under (the finitary)
Boolean operations"
\sn
\item[$(c)$]   $(p^*,r^*,q^*) \Vdash_{\bbP^{\bold x}_{\delta(*)} \times \cT_{<
\theta} \times \bbQ_{\le \theta}} ``\{\name C_i:i < \theta\} \cap
\name D$ generated $\name D$ and for some club $E$ of 
$\theta$, if $\varepsilon < \theta$ then 
$\{\name{\bold p}_\zeta:\zeta < \varepsilon\},\{\name C_i:i <
\varepsilon\} \cap \name D$ generate the same filter"
\sn
\item[$(d)$]   $E$ is actually a club of $\theta$ from $\bold V$
\sn
\item[$(e)$]  $\varepsilon \in E \Rightarrow (p^*,r^*,q^*) \Vdash
``\{\name C_i:i < \varepsilon\}$ is closed under the (finitary)
Boolean operations", so even $p^*$ forces this (for
$\Vdash_{\bbP_{\delta(*)}[\bold x_\varepsilon]}$)
\end{enumerate}
\sn
\item[$(*)_6$]   there are $r_*,q_*$ from $\cT_{< \theta},\bbQ_{\le
\theta}$ respectively and ${\cU} \in E_\theta$ such that
\begin{enumerate}
\item[$(a)$]   $\varepsilon \in {\cU} \Rightarrow
  r_\varepsilon = r_* \wedge q_\varepsilon \restriction \theta = q_*
\rest \theta$ so $r^* \le_{\cT_{< \theta}} r_\varepsilon$; also $q_* 
\le_{\bbQ_{\le \theta}} q_\varepsilon$
\sn
\item[$(b)$]   $\langle q_\varepsilon(\theta):\varepsilon \in
  {\cU}\rangle$ is a $\Delta$-system with heart $q_*(\theta) \in \bbQ_\theta$ 
\sn
\item[$(c)$]   if $\varepsilon_1 < \varepsilon_2,
\varepsilon_1 \in \cU,\varepsilon_2 \in \cU$
\then \, $q_{\varepsilon_1},q_{\varepsilon_2}$ are
compatible\footnote{so even any $< \bold f(\theta)$ members are}
\sn
\item[$(d)$]   $\cU \subseteq E$ where $E$ is from $(*)_5(d)$
\end{enumerate}
\end{enumerate}
\mn
[Why?  By the proof of Levy$(\bold f(\theta),< \theta) \models
 \theta$-c.c.]
\mn
\begin{enumerate}
\item[$(*)_7$]  for $\xi < \zeta < \theta$ let $\name D'_{\xi,\zeta}$ be
the following $\bbP^{\bold x}_{\delta(*)}$-name: 
it is the filter on $\bbN$ 
generated by the family $\{\sigma(\name C_{i_0},\dotsc,\name
C_{i_n-1}):\sigma(x_0,\dotsc,x_{n-1})$ is a Boolean term and for some
$\varepsilon \in \cU \cap \zeta \backslash \xi$ we have $\ell < n
\Rightarrow i_\ell \in (\xi,\varepsilon)$ and $\name A_\varepsilon \subseteq^*
\sigma(\name C_{i_0},\dotsc,\name C_{i_{n-1}})\}$
\sn
\item[$(*)_8$]  $\Vdash_{\bbP^{\bold x}_{\delta(*)}} ``\langle \name
D'_{\xi,\zeta}:\zeta \in (\xi,\theta]\rangle$ is increasing continuous
for each $\xi < \theta$ and $\langle \name D'_{\xi,\zeta}:
\xi < \zeta\rangle$ is decreasing for each $\zeta < \theta$
and $\emptyset \notin \name D'_{\xi,\zeta}$ for $\xi < \zeta < \theta$
and if $\xi < \zeta \in \cU$ then $\name A_\zeta,\bbN
\backslash \name A_\zeta$ are $\ne \emptyset$ mod 
$\name D'_{\xi,\zeta}"$.
\end{enumerate}
\mn
Recall $\theta < \lambda = \cf(\lambda)$ and so 
$\langle \bbP_{\delta(*)}[\bold x_\varepsilon]:\varepsilon < \lambda\rangle$ is
$\lessdot$-increasing with union $\bbP^{\bold x}_{\delta(*)}$, hence there is
$\gamma(*) < \lambda$ of cofinality $\theta$ such that for every
$\varepsilon < \theta,\bold q_\varepsilon,\bold q'_\varepsilon$ 
are representations of $\bbP_{\delta(*)}[\bold x_{\gamma(*)}]$-name 
so $\name A'_\varepsilon,\name C_\varepsilon$ are 
$\bbP_{\delta(*)}[\bold x_{\gamma(*)}]$-names and let
$\bold j_{\gamma(*)}$ be the $\bold j_*$ from \ref{cn.77}, so $(\bold
j_{\gamma(*)},\bold x_{\gamma(*)},\bold x_{\gamma(*)+1})$ here stand for
$(\bold j_*,\bold x,\bold y)$ there.

Recall $\langle \bbP_{\delta(*)}[\bold x_\varepsilon]:\varepsilon <
\lambda\rangle$ is $\lessdot$ increasing and is
continuous for ordinals of cofinality $> \aleph_0$.  
Let $\name A'_\theta$ be $\bold j_{\gamma(*)}
(\langle \name A'_\varepsilon:\varepsilon \in \cU\rangle/E_\theta)$, 
well abusing our notation a little; you may prefer to use $\bold
q_\theta = \bold j_{\gamma(*)}(\langle \bold q_\varepsilon:\varepsilon
< \theta \rangle/E_\theta)$ and $\name A'_\theta$ be the 
$\bbP_{\delta(*)}[\bold x_{\gamma(*)+1}]$-name represented by $\bold q_\theta$.

Now as $(p^*,r^*,q^*) \Vdash ``\{\name C_i:i < \theta\} \cap \name D$ 
generate an ultrafilter on $\bbN"$ and $(p^*,r^*,q^*)$ is below
$(p^*,q_{\min(\cU)},r_{\min(\cU)})$ so there is $(p^1,r^1,q^1) \in
\bbP^{\bold x}_{\alpha(*)} * \cT_{< \theta} * \bbQ_{\le \theta}$ 
above it, $n \in \bbN,\varepsilon_0,
\dotsc,\varepsilon_{n-1} < \theta$, Boolean term
$\sigma(x_0,\dotsc,x_{n-1})$ and truth value $\bold t$ such that 
\mn
\begin{enumerate}
\item[$(*)_9$]  $(p^1,r^1,q^1)$ forces $\sigma(\name
C_{\varepsilon_0},\dotsc,\name C_{\varepsilon_{n-1}}) \in \name D$ and
is included in $(\name A'_\theta)^{[\bold t]}$ recalling $A^{[1]} =
A,A^{[0]} = \bbN \backslash A$ 
\end{enumerate}
\mn
\underline{hence}
\mn
\begin{enumerate}
\item[$(*)_{10}$]  $p^1 \Vdash_{\bbP_{\delta(*)}[\bold x]} 
``\sigma(\name C_{\varepsilon_0},\dotsc,\name C_{\varepsilon_{n-1}})
\subseteq^* (\name A'_\theta)^{[\bold t]}"$.
\end{enumerate}
\mn
Let $p^2 \in \bbP_{\delta(*)}[\bold x_{\gamma(*)+1}]$ 
be such that $p^2 \le p^* \in \bbP_{\delta(*)}[\bold x_{\gamma(*)+1}] 
\Rightarrow p^1,p^*$ compatible, so clearly
\mn
\begin{enumerate}
\item[$(*)_{11}$]  $p^2 \Vdash_{\bbP_{\alpha(*)}[\bold x_{\gamma(*)+1}]} 
``\sigma(\name C_{\varepsilon_0},\dotsc,\name C_{\varepsilon_{n-1}})$ is
$\subseteq^*(\name A'_\theta)^{[\bold t]}"$.
\end{enumerate} 
\mn
Let $\langle p^2_\varepsilon:\varepsilon < \theta\rangle \in
{}^\theta(\bbP_{\delta(*)}[\bold x_{\gamma(*)}])$ be such that $\bold
j_{\gamma(*)}(\langle p^2_\varepsilon:\varepsilon < \theta\rangle) = p^2$.  

Hence
\mn
\begin{enumerate}
\item[$(*)_{12}$]   $\cU_1 = \{\zeta \in \cU:p^2_\zeta \Vdash ``\sigma(\name
C_{\varepsilon_0},\dotsc,\name C_{\varepsilon_{n-1}})$ is
$\subseteq^* (\name A'_\zeta)^{[\bold t]}"\}$ belongs to $E_\theta$.
\end{enumerate}
\mn
Without loss of generality 
$\langle p^2_\zeta:\zeta \in \cU_1\rangle$ are pairwise
compatible hence by \L os theorem for some $\zeta$
\mn
\begin{enumerate}
\item[$(*)_{13}$]  $\zeta \in \cU_1$ so $\zeta < \theta$ and 
$p^2,p^2_\zeta$ has a common upper bound $p^3 \in
\bbP_{\bold x_{\gamma(*)+1}}$, hence $p^1,p^3$ has a common upper
bound $p^4 \in \bbP^{\bold x}_{\delta(*)}$.
\end{enumerate} 
\mn
So recalling $q_\zeta$ is from $(*)_4$,
\mn
\begin{enumerate}
\item[$(*)_{14}$]  $(p^4,r_*,q_\zeta)$ forces
\begin{enumerate}
\item[$(a)$]  $\name A'_\zeta \in \name D$
\sn
\item[$(b)$]  $\sigma(\name C_{\varepsilon_0},\dotsc,\name
C_{\varepsilon_{n-1}}) \in \name D$
\sn
\item[$(c)$]  $\sigma(\name C_{\varepsilon_0},\dotsc,\name
C_{\varepsilon_{n-1}}) \subseteq^* (\name A'_\theta)^{[\bold
t]},(\name A'_\zeta)^{[\bold t]}$.
\end{enumerate}
\end{enumerate} 
\mn
Contradiction.
\end{PROOF}

\begin{claim}
\label{e.8} 
In \ref{e.4} (and \ref{cn.19}) instead of ``$E_\theta$ is
$\theta$-complete (so $\theta$ is measurable) we may require that there is
$\Theta'_2 \subseteq \Theta_2$ such that:
\mn
\begin{enumerate}
\item[$(a)$]  $(\Theta'_2,\bold f)$ are as in \ref{e.4}
\sn
\item[$(b)$]  defining $\bbQ$ we use $\Theta'_2$ if $\theta \in
\Theta'_2$ then $E_\theta$ is $\theta$-complete
\sn
\item[$(c)$]  if $\sigma \in \Theta_2 \backslash \Theta'_2$ \then \,
$\theta = \max(\Theta'_2 \cap \sigma)$ is well defined,
$[\theta,\sigma] \cap \Theta_1 = \emptyset$ and
$E_\theta$ is a uniform $\theta$-complete ultrafilter on $\sigma$ so
$\theta$ is a $\sigma$-compact cardinal.
\end{enumerate}
\end{claim}

\begin{PROOF}{\ref{e.8}}
Similar to \ref{e.4}.
\end{PROOF}

\begin{remark}
\label{e.11}
The situation is similar for any set
$\{\aleph_\alpha:\alpha \in u\}$ of successor of regular cardinals. 
\end{remark}

\begin{claim}
\label{e.30}  
In \ref{e.16} above the sufficient conditions for
``$\theta \notin \Sp_\chi$ in $\bold V^{\bbP}$" are
sufficient also for ``$(\forall \mu)(\text{\rm cf}(\mu)) = \theta
\Rightarrow \mu \notin \Sp_\chi)"$.
\end{claim}

\begin{PROOF}{\ref{e.30}}  
The same.
\end{PROOF}

So we can resolve Problem (6) from Brendle-Shelah \cite[\S8]{BnSh:642}.
\begin{conclusion}
\label{e.34}  
If GCH and $\aleph_1 \le \theta < \kappa =
\cf(\kappa) < \lambda = \lambda^\kappa,\kappa$ is measurable,
\then \, there is a forcing notion $\Bbb P$ of cardinality $\lambda$
collapsing the cardinals in $(\theta,\kappa)$ but no others
such that in $\bold V^{\bbP}$, for every cardinal 
$\mu \in (\kappa,\lambda)$ of cofinality
$\kappa$, we have $\mu \notin \Sp_\chi \wedge \mu =
\sup(\Sp_\chi \cap \mu)$.
\end{conclusion}
\newpage

%\bibliographystyle{alphacolon}
%\bibliography{lista,listb,listx,listf,liste,listz}

\end{document}